\documentclass[10pt,english]{elsarticle}
\usepackage[T1]{fontenc}
\usepackage[utf8]{inputenc}
\usepackage{geometry}
\usepackage{color}
\usepackage{float}
\usepackage{textcomp}
\usepackage{amsmath}
\numberwithin{equation}{section}
\usepackage{amssymb}
\usepackage{graphicx}
\usepackage{setspace}
\usepackage{stmaryrd}
\usepackage{mathabx}
\usepackage{bm}
\usepackage{babel}
\usepackage{hyperref}
\usepackage{accents}
\def\n n{{n_{\rm nod}}}
\def\n p{{n_{\rm pts}}}

\newcommand{\jump}[1]{\ensuremath{[\![#1]\!]} }

\makeatletter
\def\ps@pprintTitle{%
    \let\@oddhead\@empty
    \let\@evenhead\@empty
    \def\@oddfoot{\footnotesize\itshape Preprint submitted to Elsevier\hfill November 5, 2026}%
    \let\@evenfoot\@oddfoot
    }
\makeatother

\allowdisplaybreaks

\makeatletter

\floatstyle{ruled}
\newfloat{algorithm}{tbp}{loa}
\floatname{algorithm}{\protect\algorithmname}

\makeatother

\begin{document}
\title{A posteriori space error estimates for the finite element approximation of the incompressible Navier-Stokes equations
based on the variational multiscale framework}
\author{
J. Sekhavati$\mathrm{^{a,b}}$,
I. Martínez-Suárez$\mathrm{^{a,b}}$,
S.A. Khan$\mathrm{^{a,b}}$,
J. Baiges$\mathrm{^{a,b}}$, 
G. Hauke$\mathrm{^{c}}$, 
and
R. Codina$\mathrm{^{a,b}}$
\\~\\
$\mathrm{^{a}}$ Universitat Polit\`{e}cnica de Catalunya, Jordi Girona $1-3$, Edifici C1, Barcelona 08034, Spain.
\\
$\mathrm{^{b}}$ International Center for Numerical Methods in Engineering (CIMNE), Spain.
\\
$\mathrm{^{c}}$ 
Universidad de Zaragoza, Zaragoza, 50018, Spain
\\~\\
javad.sekhavati@upc.edu (J. Sekhavati), 
ignacio.martinez.suarez@upc.edu (I. Martínez-Suárez), 
sheraz.ahmed@upc.edu (S.A. Khan), 
joan.baiges@upc.edu (J. Baiges), 
ghauke@unizar.es (G. Hauke),
ramon.codina@upc.edu (R. Codina)
}

\begin{abstract}
    This study introduces a novel a posteriori space error estimator derived within the variational multiscale (VMS) framework for the stationary incompressible Navier-Stokes equations. In this framework, the continuous problem is decomposed into coarse and fine scales. The fine scale, often referred to as the sub-grid scale (SGS), is modeled as proportional to the component of the residual orthogonal to the finite element (FE) space. This modeled SGS then serves as a reliable error estimator. Specifically, the approach employs a scaled norm of the SGSs as an a posteriori error estimate, evaluated within the stabilized norm associated with the problem. This methodology was previously applied to the stationary convection-diffusion-reaction equation \cite{codina2025posteriorierrorestimatesfinite} and is now extended to the stationary incompressible Navier-Stokes equations. The accuracy and robustness of the proposed space error estimator are assessed through a series of numerical experiments on both manufactured and benchmark incompressible flow problems. Additionally, its performance is compared in selected cases with the space error estimator proposed in \cite{Irisarri2021}.
\end{abstract}

\begin{keyword}
    \texttt{A posteriori space error estimator, Variational multiscale method, Incompressible Navier-Stokes equations}
\end{keyword}

\maketitle

\section{Introduction}
\noindent
Modeling incompressible fluid flows through the Navier-Stokes equations is crucial in various scientific and technical fields, such as aerodynamics, biomedical flows, and weather forecasting. Approximating these flows numerically presents several computational challenges of various natures. One of the most relevant challenges in complex flow problems is balancing accuracy and computational efficiency to obtain reliable numerical simulations. This work focuses on introducing an a posteriori space error estimator with the objective of guiding an adaptive finite element (FE) simulation strategy for fluid mechanics problems.

Employing adaptive mesh strategies guided by a posteriori error estimation has been an active area of research in computational mechanics \cite{ainsworth1997posteriori}. The a posteriori error estimator in space assesses the spatial distribution of numerical errors in the computational domain and serves as an essential tool for adaptive mesh refinement (AMR) algorithms \cite{bangerth2003adaptive,verfurth1994posteriori}. The core idea is to optimize computational efficiency by concentrating the fine discretization and, thus, the computational effort in complicated flow regions that give high error estimates, while the smoother regions are modeled with coarse discretizations.

Traditional AMR is based on simple approaches such as geometric criteria, physics-based gradient criteria, or Zienkiewicz-Zhu (ZZ)~\cite{zienkiewicz1987simple}. However, advanced error-based criteria can lead to more efficient refinement. In the context of the Navier-Stokes and related fluid equations, residual-based a posteriori error estimators have been extensively investigated \cite{bernardi2009finite,bernardi2016posteriori,jin1998posteriori,prudhomme1999posteriori,ervin1996posteriori,burda2007numerical,williamson2019posteriori}. For instance, Burda et al. \cite{burda2007numerical} examined residual-based error estimators to analyze accuracy in stabilized finite element formulations of flow problems, while Williamson et al. \cite{williamson2019posteriori} proposed residual-based a posteriori estimates for driving adaptive mesh refinement in coupled flow systems like the Stokes-Brinkman problem. In \cite{offermans2020adaptive}, Offermans et al. developed AMR strategies driven by goal-oriented error estimators based on a functional quantity of interest and compared their performance against local spectral error estimators for the Navier-Stokes equations. Moreover, Berrone and Marro \cite{berrone2009space,berrone2010numerical} established a space and time adaptive strategy for time-dependent Navier--Stokes flows. By independently employing a residual-based error indicator for space and a velocity-gradient error indicator for time, they demonstrated reliable error control and stable effectivity indices across their selected benchmark problems. Recent developments in adaptive mesh refinement for incompressible flow simulations have also explored alternative discretization frameworks, such as THB-spline-based isogeometric analysis \cite{bastl2025adaptive}.

In recent years, a posteriori space error estimation approaches derived from the VMS method have been extensively applied to a wide range of computational problems. The VMS approach established by Hughes and his collaborators \cite{hughes1995multiscale,hughes1998variational} has greatly contributed to the development of stabilized FE methods. The key point of the VMS approach is to split the unknowns of the problem (velocity and pressure) into their FE components and a sub-grid scale (SGS) component, which, in fact, needs to be approximated.

The VMS methodology has been applied across various scientific disciplines, including the development of stabilized FE formulations \cite{guermond1999stabilization,masud2004multiscale,araya2006stabilized}, large-eddy simulation turbulence modeling \cite{codina2018variational,bazilevs2007variational,ahmed2017review}, and a posteriori error estimation techniques. The a posteriori error estimation approach based on the VMS framework models the SGS as an estimator of the numerical error. This concept was introduced in the early development of the VMS framework \cite{hughes1998variational}, and the first attempt to derive a posteriori error estimates within this theory was presented in \cite{hauke2006multiscale} for convection-dominated problems. The same strategy is applied to derive the VMS-based a posteriori error estimators for the Navier-Stokes equations in \cite{hauke2015variational, Irisarri2021}. In \cite{rossi2013parallel}, Rossi et al. also used a VMS-based error estimator to implement a parallel AMR algorithm for incompressible flow simulations. In this work, the VMS-based error estimation framework, previously applied to reduced order models \cite{codina2021posteriori}, particular cases of the compressible Navier-Stokes equations \cite{bayona2018variational}, solid mechanics problems \cite{baiges2017variational}, and convection-diffusion-reaction problems \cite{khan2025variationalmultiscaleapproachgoaloriented,codina2025posteriorierrorestimatesfinite}, is extended to the incompressible Navier-Stokes equations. A comprehensive review of the development of VMS-based a posteriori error estimation techniques in fluid mechanics is presented in \cite{hauke2023review}. 

In this study, we present an a posteriori space error estimator for the incompressible Navier-Stokes equations based on the VMS methodology. The error estimate is derived using the orthogonal SGS (OSGS) method introduced in \cite{codina2000stabilization}, where the SGSs are chosen to be orthogonal to the FE space. For space error control, the SGS model is used to construct reliable estimators. The VMS-based stabilized FE formulation is solved using a monolithic approach.

Another contribution of this work is that we compare our approach, which uses the same velocity and pressure SGSs as for the FE calculation, with the approach proposed in \cite{Irisarri2021}. There, given a FE solution, a velocity SGS is reconstructed, and its $H^1$ seminorm subsequently serves as an a posteriori error estimator. In contrast, the a posteriori error estimator we propose is only meaningful when the error is measured in the norm associated with the VMS-stabilized FE formulation, often simply called a stabilized norm.

The remainder of the paper is organized as follows. In Section 2, the incompressible Navier-Stokes problem is stated, and its FE formulation is summarized, along with the VMS-based stabilization technique. Section 3 introduces our a posteriori space error estimator in the stabilized norm, together with related $L^2$ and $H^1$ estimates that are evaluated for comparison. In Section 4, the reference space error estimator in the $H^1$ seminorm from Irisarri and Hauke \cite{Irisarri2021} is briefly presented. Section 5 presents the numerical examples to illustrate the performance
of the stabilized FE formulation and the proposed a posteriori error estimator. Moreover, selected cases are included to compare the behavior of the proposed space error estimator with the reference one. Finally, the conclusions of this study are presented in Section 6.   

\section{Finite element approximation of the incompressible Navier-Stokes equations using the variational multiscale method}

\subsection{Flow equations}
\noindent In this section, we summarize the VMS-OSGS approach applied to the incompressible Navier-Stokes equations described in \cite{codina2002stabilized,codina2007time}. The incompressible Navier-Stokes equations consist of finding a velocity $\bm{u}:\Omega \longrightarrow\mathbb{R}^d$, and a pressure $p:\Omega  \longrightarrow \mathbb{R}$ such that
\begin{eqnarray}
-2\nabla\cdot\mu\nabla^{S}\bm{u}+\rho \bm{u}\cdot\nabla\bm{u}+\nabla p & =\bm{f} & \quad\mbox{ in }~\Omega,\label{eq:strns1}\\
\nabla\cdot\bm{u} & =0 & \quad\mbox{ in }~\Omega,\\
\bm{u} & =\bar{\bm{u}} & \quad\mbox{ on }~\Gamma_{D},\\
\bm{n}\cdot\bm{\sigma} & =\bm{0} & \quad\mbox{ on }~\Gamma_{N}.\label{eq:strns}
\end{eqnarray}
where $\nabla^{S}\bm{u}=\frac{1}{2}(\left(\nabla\bm{u}\right)^{T}+\nabla\bm{u})$ and ${\bm\sigma} = -p {\bm I} + 2\mu\nabla^S \bm{u}$ is the 
Cauchy stress tensor.
$\Omega\subset\mathbb{R}^{d}$ is a bounded domain of dimension $d$, $\mu$ is the
dynamic viscosity, $\rho$ is the density and $\bm{f}$ is a body force. $\Gamma_{D}$ is the
Dirichlet boundary, where velocity boundary conditions are applied, and $\Gamma_{N}$ is the Neumann boundary, where the traction is prescribed to zero (for simplicity) to represent a traction-free outflow boundary.

Let us now consider the variational weak form of equations
(\ref{eq:strns1})-(\ref{eq:strns}). We define $V=H^{1}(\Omega)^{d}$, $V_{0}=\{\bm{v}\in V\mbox{ }|\mbox{ }\bm{v}=\bm{0}$ on $~\Gamma_{D}\}$, and $Q=L^{2}(\Omega)$. In the following, $(\cdot,\cdot)$ stands for the $L^{2}(\Omega)$ inner product and $\left\langle \cdot,\cdot\right\rangle $ for the integral of
the product of two functions, not necessarily in $L^{2}(\Omega)$. The variational problem consists of finding
$[\bm{u},p]\in V\times Q$ such that:
\begin{align}
\langle\bm{v},\rho\bm{u}\cdot\nabla\bm{u}\rangle+2\mu\left(\nabla^{S}\bm{v},\nabla^{S}\bm{u}\right)-\left(p,\nabla\cdot\bm{v}\right)=\left\langle \bm{v},\bm{f}\right\rangle \quad\forall\bm{v}\in V_{0},\label{eq:wkns1}\\
\left(q,\nabla\cdot\bm{u}\right)=0\quad\forall q\in Q,\label{eq:wkns2}
\end{align}
with
\begin{align*}
\bm{u}=\bar{\bm{u}}\quad\mbox{ on }~\Gamma_{D}.
\end{align*}

Let us define
\[
B\left([\bm{v},q],[\bm{a};\bm{u},p]\right):=\langle\bm{v},\rho\bm{a}\cdot\nabla\bm{u}\rangle+2\mu\left(\nabla^{S}\bm{v},\nabla^{S}\bm{u}\right)-\left(p,\nabla\cdot\bm{v}\right)+\left(q,\nabla\cdot\bm{u}\right),
\]
where $\bm{a}$ represents the convective velocity. Given $\bm{a}$, $B$ is a bilinear form on $(V\times Q)\times (V\times Q)$. The problem to be solved corresponds to taking $\bm{a} = \bm{u}$, so that $B\left([\bm{v},q],[\bm{u};\bm{u},p]\right)$ is semilinear. It is easily verified that it is continuous.

\subsection{Finite element approximation}

\noindent In the Galerkin FE framework, the semilinear form 
$B$ is known to exhibit stability issues arising from the convective nature of the flow when the element Reynolds number is high, and additionally it has compatibility requirements between the velocity and pressure spaces derived from the classical inf-sup condition.
In this work, we use the VMS method to stabilize the FE approximation, but also to model the approximation error by means of the SGSs, which we describe next.

Let us start by considering the FE partition $\mathcal{T}_{h}:=\{K\}$ defined over the
domain $\Omega$, $h_K$ denoting the element size, $h = \max_K h_K$ and $N$ being the total number of elements. From this
partition we construct the FE spaces $V_{h}\subset V,$
$Q_{h}\subset Q$. The VMS method is defined by introducing the SGS space for the velocity field $\tilde{V}$, and the pressure SGS space $\tilde{Q}$ such that:
\begin{align*}
V =  V_{h}\oplus\tilde{V}, \quad Q  =  Q_{h}\oplus\tilde{Q}.
\end{align*}
This implies that the velocity and pressure solution fields are now decomposed into the FE approximation, plus a model for the part of the solution which cannot be captured by the FE space, the SGSs. The contribution of the SGSs will serve to stabilize the formulation and enhance the obtained solution, and in the present work also as an error estimator for the spatial accuracy of the numerical approximation, as we shall see. The splitting of the unknowns we consider is therefore:
\begin{align}
\bm{u} = \bm{u}_{h}+\tilde{\bm{u}}, \quad p = p_{h}+\tilde{p}.\label{eq:split}
\end{align}

The SGSs $\tilde{\bm{u}}$ and $\tilde{p}$ must now be approximated. In the case of the velocity SGS $\tilde{\bm{u}}$, an approximation will be used when it is evaluated in the element interiors, which will be denoted $\tilde{\bm{u}}_K$, and another approximation will be used when it is evaluated on the edges $E$ of the FE partition, which will be denoted $\tilde{\bm{u}}_E$. Thus, using the splitting \eqref{eq:split} in \eqref{eq:wkns1}-\eqref{eq:wkns2} and integrating by parts, one finds that the original problem is equivalent to finding $[\bm{u}_{h},p_{h}]\in V_{h}\times Q_{h}$, as well as the SGSs $\tilde{\bm{u}}_K$, $\tilde{\bm{u}}_E$ and $\tilde{p}$, such that
\begin{eqnarray}
B\left([\bm{v}_{h},q_{h}],[\bm{u};\bm{u}_{h},p_{h}]\right)
+\sum_{K}\left\langle -\mu\Delta\bm{v}_{h}-\rho\bm{u}\cdot\nabla\bm{v}_{h}-\nabla q_{h},\tilde{\bm{u}}_K\right\rangle _{K},\label{eq:u_h_equation_1}\nonumber\\
+\sum_{E}\left\langle \jump{2\mu\bm{n}\cdot\nabla^S\bm{v}_{h}+q_h\bm{n}} ,\tilde{\bm{u}}_{E}\right\rangle _{E} \nonumber\\
-\sum_{K}\left\langle \nabla\cdot\bm{v}_{h},\tilde{p}\right\rangle _{K} & = & \left\langle \bm{v}_{h},\bm{f}\right\rangle ,\label{eq:u_h_equation_3}
\end{eqnarray}
for all test functions $\bm{v}_{h},$ $q_{h}$ and where $\left\langle \cdot,\text{\ensuremath{\cdot}}\right\rangle _{K}$ stands for the integral of the product of two functions in $K$, and $\left\langle \cdot,\text{\ensuremath{\cdot}}\right\rangle _{E}$ for the integral of the product of two functions over a face (edge in 2D) $E$. The jump operator across a face $E$ shared by two elements $K_1$ and $K_2$ is defined as: 
\begin{eqnarray}
\jump{{\bm{n}}g}_{E}=\bm{n}_{1}g\vert_{\partial K_{1}\cap E}+\bm{n}_{2}g\vert_{\partial K_{2}\cap E}, \\
\jump{\bm n\cdot \nabla g}_E = \bm{n}_{1} \cdot\nabla g\vert_{\partial K_1\cap E} + \bm{n}_{2}\cdot\nabla g\vert_{\partial K_2\cap E},
\end{eqnarray}
where $\bm{n}_{i}$ denotes the outward unit normal on the boundary of element $K_{i}$.

Note that in \eqref{eq:u_h_equation_3} the advection velocity is $\bm{u}$, which can be computed as $\bm{u}_h + \tilde{\bm{u}}_K$ using the approximation of the velocity SGS given below. However, it can also be approximated as $\bm{u}_h$. See, e.g., \cite{colomes2015assessment} for a discussion about the implications of one choice or another. 

The approximation of the SGSs needs to be introduced now. Here, it is important to highlight that the SGSs can be used for both the FE calculation and for the a posteriori error estimation. In the first a posteriori error approach to be discussed, {\em they are the same}, but in Section~\ref{sec:sgs-recov} we shall present a method in which the SGSs are recovered from the FE solution and used to estimate the error.

Let us start with the SGSs approximation used both in FE calculations and in the a posteriori error estimation. We shall not motivate this approximation, but rather write the final expressions of the approximated SGSs. Different ways to motivate them can be found in the review \cite{codina2017VariationalMultiscale}. The expression of the approximation we use for the SGSs, without using any particular notation for this approximation, is
\begin{align}
\tilde{\bm{u}}_K & =\tau_{u,K} P_{\tilde{V}}\left[\bm{R}_h\right]_K 
=\tau_{u,K} P_{\tilde{V}}\left[\bm{f}-\left(-\mu\Delta\bm{u}_{h}+\rho\bm{u}\cdot\nabla\bm{u}_{h}+\nabla p_{h}\right)\right]_K,\label{eq:navier-stokes_subscales} \\
\tilde{\bm{u}}_{E} & =-\tau_{E}\jump{\bm{n}\cdot\bm{\sigma}_h}_E=-\tau_{E}\jump{2\mu\bm{n}\cdot\nabla^S\bm{u}_h-p_{h}\bm{n}}_{E},\label{eq:boundary_subscales} \\
\tilde{p}_K & = - \tau_{p,K}P_{\tilde{Q}}\left(\nabla\cdot\bm{u}_{h}\right).\label{eq:pressure_subscales}
\end{align}

Let us explain the different terms appearing in these expressions. $\left[\bm{R}_h\right]_K$ is the residual of the momentum equation evaluated in element $K$, $-\nabla\cdot\bm{u}_{h}$ can be understood as the residual of the continuity equation, and $\jump{\bm{n}\cdot\bm{\sigma}_h}_E$ is the jump of the normal component of the stress on edge $E$. Only interior edges and edges on $\Gamma_N$ need to be considered; for the latter, $\jump{2\mu\bm{n}\cdot\nabla^S\bm{u}_h-p_{h}\bm{n}}_{E} = {2\mu\bm{n}_E\cdot\nabla^S\bm{u}_h \vert_E-p_{h}\vert_E\bm{n}_E}$, as the prescribed traction is zero. Thus, the approximation presented can be considered as residual-based. 

$P_{\tilde{V}}$ and $P_{\tilde{Q}}$ are the projections onto the spaces of SGSs $\tilde{V}$ and $\tilde{Q}$, respectively, that need to be chosen. Let $P_h$ be the $L^2$ projection onto the FE space, either $V_h$ when the argument is a vector or onto $Q_h$ when the argument is a scalar. We will consider two different choices, yielding two different VMS formulations. In the case of the OSGS (see \cite{codina2000stabilization}), we take 
\begin{eqnarray}
\tilde{V}=V_{h}^{\perp}\cap V, \qquad P_{\tilde{V}}=I-P_h =: P_h^\bot \quad \hbox{(applied to vectors)},\\
\tilde{Q}=Q_{h}^{\perp}\cap Q, \qquad P_{\tilde{Q}}=I-P_h = P_h^\bot \quad \hbox{(applied to scalars)}.
\end{eqnarray}
The other option we shall consider is what we call the Algebraic SGS (ASGS) formulation, in which
\begin{eqnarray}
\tilde{V}= \hbox{space of FE residuals of $V_h$}, \qquad P_{\tilde{V}}=I\quad \hbox{(applied to vectors)},\\
\tilde{Q}= \hbox{space of FE residuals of $Q_h$}, \qquad P_{\tilde{Q}}=I\quad \hbox{(applied to scalars)}.
\end{eqnarray}

It only remains to define the stabilization parameters $\tau_{u,K}$, $\tau_E$ and $\tau_{p,K}$ in \eqref{eq:navier-stokes_subscales}-\eqref{eq:pressure_subscales}. If $h_K$ is the characteristic length of element $K$ and $h_E$ the characteristic length of edge $E$, we compute them as
\begin{align*}
\tau_{u,K}^{-1} & =\left(c_{1}\frac{\mu}{h_K^{2}}+c_{2}\frac{\rho|\bm{u}|}{h_K}\right),\\
\tau_{p,K} & =c_{3}h_K^{2}\tau_{u,K}^{-1},  \\
\tau_{E} & =\frac{\delta_{0}h_E}{\mu + \rho|\bm{u}| h_E},
\end{align*}
where $c_1$, $c_2$, $c_3$ and $\delta_0$ are dimensionless algorithmic constants. For linear elements, we take them as  $c_1 = 4$, $c_2 = 2$, $c_3 = 1$ and $\delta_0 = 1/12$. If $k$ is the polynomial order of the FE interpolation, we scale them with $k^4$, $k^2$, $k^{-4}$ and $k^{-2}$, respectively. The velocity $|\bm{u}|$ is a characteristic value, either in $K$ or in $E$, for example, the maximum. 

This completes the required expression for the SGSs in terms of the FE unknowns. Note that (\ref{eq:navier-stokes_subscales})-(\ref{eq:pressure_subscales}) need to be evaluated at the numerical integration points, since the SGSs are needed there to compute (numerically) the integrals appearing in (\ref{eq:u_h_equation_1}). Let us point out that, while these are the expressions of the SGSs to be used in the a posteriori error analysis (see below), sometimes they are not needed in the FE calculation. Indeed, in many cases one can set $\tau_{p,K} = 0$, and $\tau_E$ is only needed when the pressure interpolation is discontinuous. For stability reasons, the only fundamental term is $\tau_{u,K}$.

\section{Proposed space error estimator in the stabilized norm}

\noindent Let us introduce our first proposal to obtain an a posteriori error estimate using the SGSs associated with the VMS formulation. This first approach states that {\em a properly scaled norm of the SGSs provides an a posteriori error estimate in the stabilized norm of the formulation.} This idea was initiated in \cite{baiges2017variational,codina2021posteriori} and elaborated in \cite{codina2025posteriorierrorestimatesfinite}. A similar approach for the convection-diffusion equation approximated using the SUPG formulation was analyzed in \cite{john-novo-2013}. The analysis in \cite{codina2025posteriorierrorestimatesfinite} shows that, under certain assumptions, the resulting estimate is indeed robust, i.e., it bounds from above and from below the true error.

We have to define the stabilized norm of the formulation and the scaled norm of the SGSs. Let us start with the former. The stabilized norm is the norm in which it is possible to prove stability and convergence. For the incompressible Navier-Stokes equations, it is given by \cite{codina2000stabilization}:
\begin{align}
    |||[\bm{v},q]|||^2=\mu\|\nabla \bm{v}\|^2+\sum_K\tau_{u,K}\| \rho \bm{a}\cdot\nabla\bm{v} +\nabla q\|^2_K+\sum_K\tau_{p,K}\|\nabla \cdot \bm{v}\|^2_K + \sum_E  \tau_E \Vert \jump{2\mu\bm{n}\cdot\nabla^S\bm{v}-q\bm{n}}\Vert_E^2.
    \label{eq:stabilized_norm}
\end{align}
When applied to a FE solution $[\bm{u}_h, p_h]$, the advection velocity $\bm{a}$ can be either $\bm{u}$ or $\bm{u}_h$. Let us point out that the last term, $\sum_E  \tau_E \Vert 2\mu\bm{n}\cdot\nabla^S\bm{v}-q\bm{n}\Vert_E^2$, can be bounded from above by the first two terms using a standard trace inequality, so that in practice it is not needed, and therefore not used in the numerical examples presented here.

The stabilized norm of the exact error $\bm{E} = [\bm{u}-\bm{u}_h,p-p_h]$ of the FE solution is
\begin{align}
    |||\bm{E}|||=|||[\bm{u}-\bm{u}_h,p-p_h]|||\,.
\end{align}

Let us now present the a posteriori error estimate that we propose. Assuming for simplicity constant stabilization parameters within each element $K$, the element contribution to the a posteriori error that we propose is:
\begin{equation}
\eta_{K}^{2}( \tilde{\bm{u}}_K,\tilde{p}_K,\tilde{\bm{u}}_{E}):=\tau_{u,K}^{-1}\int_{K}\vert\tilde{\bm{u}}_{K}\vert^{2}+\tau_{p,K}^{-1}\int_{K}\vert\tilde{p}_K \vert^{2}+ \sum_{E\subset \partial K}\tau_{E}^{-1}\int_{E}\vert\tilde{\bm{u}}_{E}\vert^{2}\label{eq:eq4}
\end{equation}

To compute the edge term of \eqref{eq:eq4}, we do an additional approximation. Following \cite{El-Alaoui-2004}, and as done in \cite{badia2013adaptive,baiges2017variational}, the SGSs on the edges term can be shown to be equivalent to
$\Vert\bm{\sigma}_{h}-P_{h}(\bm{\sigma}_{h})\Vert_{K}^{2}$ in the sense that
\begin{equation}
\gamma_{1} \sum_{E\subset \partial K} h_{K}\int_{E}\vert\jump{\bm{\sigma}_h\cdot\bm{n}}_{E}\vert^{2}\leq\Vert\bm{\sigma}_{h}-P_{h}(\bm{\sigma}_{h})\Vert^{2}_K\leq\gamma_{2}\sum_{E\subset \partial K} h_{K}\int_{E}\vert\jump{\bm{\sigma}_h\cdot\bm{n}}_{E}\vert^{2}\label{eq:boundary_subscales_equivalency}
\end{equation}
for certain positive constants $\gamma_{1}$, $\gamma_{2}$. This allows us to approximate: 
\begin{align}
\sum_{E\subset \partial K}
\tau_{E}^{-1}\int_{E}\vert\tilde{\bm{u}}_{E}\vert^{2}\approx\tau_{E}\gamma h_{K}^{-1}\Vert{P_h^\bot}(\bm{\sigma}_{h})\Vert_{K}^{2},\label{eq:equiv}
\end{align}
for a constant $\gamma$ that needs to be determined from numerical experiments (taken as $\gamma = 18$ in our numerical examples). 

Finally, the total estimated space error in the computational domain $\eta_{\rm stab}$ will be
\begin{align}
\eta_{\rm stab}^2=\sum_K\eta_{K}^{2}( \tilde{\bm{u}}_K,\tilde{p}_K,\tilde{\bm{u}}_{E}).\label{eq:eta-stab}
\end{align}

Our claim is that this is a robust a posteriori space error estimate, which means that there exist positive constants $E_1$ and $E_2$ such that
\begin{align}
E_1 \eta_{\rm stab} \leq |||\bm{E}||| \leq E_2 \eta_{\rm stab}.\label{eq:main-claim}
\end{align}
The sharpness of these bounds can be measured by the effectivity index of the error estimator, $\mathcal{I}_{\rm eff}$, defined as
\begin{align}
\mathcal{I}_{\rm eff}= \frac{\eta_{\rm stab} }{|||\bm{E}|||} = 
\frac{\eta_{\rm stab} }{|||[\bm{u}-\bm{u}_h,p-p_{h}]|||}. \label{eq:eff-i}
\end{align}
If the claim \eqref{eq:main-claim} is correct, $\mathcal{I}_{\rm eff}$ should tend to a positive constant as $h\to 0$. In the ideal case, this constant should be 1.

As noted above, our claim is that the SGSs provide an a posteriori error estimate only in the stabilized norm. However, we could wonder whether it could be possible to obtain also a posteriori error estimates in $L^2$ and $H^1$ from the SGSs.
Later in the numerical example section, the feasibility of the presented SGSs approximation as a space error estimator will be evaluated not only in the stabilized norm, but also in the $L^2$ norm and $H^1$ seminorm of the velocity error (not the pressure), so they are worth defining. For the $L^2$ norm, a straightforward velocity space error estimate and exact error can be defined as
\begin{align}
\eta_{L^2} & =\left(\sum_K \Vert\tilde{\bm{u}}_{K}\Vert^2_{L^2(K)}\right)^{1/2},\label{eq:eta-L2}\\
E_{L^2} & =\Vert\bm{u}-\bm{u}_{h}\Vert_{L^2(\Omega)},\nonumber
\end{align}
respectively, and for the $H^1$ seminorm, the space error estimator and exact error can be defined respectively as
\begin{align}
\eta_{H^1} & = \left(\sum_K\|\nabla \tilde{\bm{u}}_{r,K}\|^2_{L^2(K)}\right)^{1/2},\label{eq:sgs-eta-h1}\\
E_{H^1} & =|\bm{u}-\bm{u}_{h}|_{H^1(\Omega)}.\label{eq:e-H1}
\end{align}
Since the approximated SGSs $\tilde{\bm u}_K$ are in general discontinuous between elements, we have introduced $\bm{u}_{r,K}$ in \eqref{eq:sgs-eta-h1}, which is a smoothed version of $\tilde{\bm{u}}_K$. We compute this smoothed velocity SGS as a global $L^2$ projection of the element SGSs $\tilde{\bm{u}}_K$ onto $V_h$, from which the gradient can be computed and localized to each element $K$.

\section{An a posteriori error estimate in the $H^1$ seminorm based on a reconstruction of the sub-grid scales}\label{sec:sgs-recov}

\noindent In this section we present an alternative way to compute an a posteriori error estimate in space for the incompressible Navier-Stokes equations specifically designed to be computed in the $H^1$ velocity seminorm. Again, the starting point is that the velocity SGS has to be proportional to the residual, but instead of computing it pointwise, its components onto the space of SGS are computed, and the scaling factors are properly designed. This space is assumed to be spanned by bubble functions within each element. This error estimate is presented in the paper by Irisarri and Hauke~\cite{Irisarri2021} (see also references therein).

To present the estimate, let us consider the 2D case, with coordinates $x$ and $y$. Within each element $K$, let us define the residuals:
\begin{align}
\mathcal{R}_{Mx} &= f_{x}-\rho \bm{u}_h\cdot\nabla u_{x,h} -\partial_{x}p_{h}+\mu\Delta u_{x,h} \quad \text{on} ~ K, \label{eq:RMx}\\
\mathcal{R}_{My} &= f_{y}-\rho \bm{u}_h\cdot\nabla u_{y,h} -\partial_{y}p_{h}+\mu\Delta u_{y,h} \quad \text{on} ~ K, \label{eq:RMy}\\
\mathcal{R}_{C} &= -\nabla\cdot \bm{u}_{h} \quad \text{on}~ K,\\
\mathcal{R}_{S_{\perp}} &= \mu\jump{\nabla\bm{u}_{h}\cdot\bm{n}}\cdot \bm{n}
\quad \text{on} ~\partial K,\\
\mathcal{R}_{S_{\parallel}} &= \mu\jump{\nabla\bm{u}_{h}\cdot\bm{n}}\cdot \bm{n}_{\parallel}
\quad \text{on} ~ \partial K.
\label{eq:RStangent}
\end{align}
These are respectively the residuals of the momentum equation along the $x$ and $y$ directions, the residual of the continuity equation, the normal component of the jump of the normal velocity gradient on the element boundary and the tangent one. The unit vector $\bm{n}_{\parallel}$ is tangent to $\partial K$. 

Consider a fixed element $K$ of the FE partition. Let $E_j$, $j = 1,\dots,n_{EK}$ be the edges of element $K$, so that $\partial K = \bigcup E_j$. Let $\ell_j$ be the length of $E_j$ and $\vert K\vert$ the area of $K$. On $K$ let us define a set of bubble functions $v_{{\rm bub},i}$, $i = 1,\dots, {n_{\text{bub}}}$, and on $E_j$ a set of edge bubble functions, still denoted $v_{{\rm bub},i}$, $i = 1,\dots, {n_{\text{edge}}}$. The idea of the a posteriori error estimator proposed in \cite{Irisarri2021} is to reconstruct the velocity SGS by considering it proportional to the residuals defined above and spanned by the bubble functions described. Thus, the velocity SGS is taken as the projection of the residuals onto the bubbles and its $H^1$ seminorm is computed as:
\begin{equation}
\begin{aligned}
|\tilde{\bm{u}}_{\rm{rec}}|_{H^1(K)} := {} &
\frac{1}{\sqrt{|K|}}
(\tau)_{H^1}^{{\mathcal{R}_{Mx}}}
\sum_{i=1}^{n_{\text{bub}}}
\left|(\mathcal{R}_{Mx}, v_{{\rm bub},i})_{K}\right| \\
&+ \frac{1}{\sqrt{|K|}}
(\tau)_{H^1}^{{\mathcal{R}_{My}}}
\sum_{i=1}^{n_{\text{bub}}}
\left|(\mathcal{R}_{My}, v_{{\rm bub},i})_{K}\right| \\
&+ \frac{1}{\sqrt{|K|}}
(\tau)_{H^1}^{{\mathcal{R}_{C}}}
\sum_{i=1}^{n_{\text{bub}}}
\left|(\mathcal{R}_{C}, v_{{\rm bub},i})_{K}\right| \\
&+ \sqrt{|K|}
(\tau)_{H^1}^{{\mathcal{R}_{S_\perp}}}
\sum_{j=1}^{n_{EK}}\sum_{i=1}^{n_{\text{edge}}}
\frac{1}{\ell_j}
\left|(\mathcal{R}_{S_\perp}, v_{{\rm bub},i})_{E_j}\right| \\
&+ \sqrt{|K|}
(\tau)_{H^1}^{{\mathcal{R}_{S_\parallel}}}
\sum_{j=1}^{n_{EK}}\sum_{i=1}^{n_{\text{edge}}}
\frac{1}{\ell_j}
\left|(\mathcal{R}_{S_\parallel}, v_{{\rm bub},i})_{E_j}\right|.
\end{aligned}
\label{eq:GH_Space_Error_Indicator}
\end{equation}

In these expressions, it remains to define how to compute the scaling factors $(\tau)_{H^1}^{\mathcal{R}_{i}}$, $i = Mx, My, C, S_{\perp}$ and $S_{\parallel}$. They can be interpreted as the mean of the gradient of the velocity SGS {\em subject to a unit residual $\mathcal{R}_{i}$} for each $i$ and are evaluated as
\begin{align}
    (\tau)_{H^1}^{\mathcal{R}_{i}}=\min\left\{ (\tau)_{\rm conv}^{\mathcal{R}_{i}} , (\tau)_{\rm visc}^{\mathcal{R}_{i}}\right\},
    \label{eq:time-scales}
\end{align}
where $(\tau)_{\rm conv}^{\mathcal{R}_{i}}$ is the expression needed in the convective limit and $(\tau)_{\rm visc}^{\mathcal{R}_{i}}$ in the viscous-dominated limit. 

The convective limit is estimated simply by solving a transport equation with a unit forcing term. The gradient of the velocity SGS must scale as $1/U$, where $U$ is the magnitude of the transport velocity. In our case, we take it as $\rho \vert \bm{u}_h\vert$, where $\vert \bm{u}_h\vert$ is a characteristic value of the velocity norm in $K$. To match the approximation with the exact SGS in 1D, we take
\begin{align}
(\tau)_{\rm conv}^{\mathcal{R}_{i}} = 
\begin{cases}\frac{1}{\rho |\bm{u}_{h}|}\frac{\alpha}{\sqrt{3}} & \hbox{for}~ i = Mx, My\\
\infty & \hbox{for}~ i = C, S_\perp, S_\parallel 
\end{cases}\label{eq:tau-conv} 
\end{align}
where $\alpha=h_K\rho|\bm{u}_{h}|/(2\mu)$ is the local Reynolds number. 

To compute $(\tau)_{\rm visc}^{\mathcal{R}_{i}}$, let us consider {\em a Stokes problem posed in the reference element $\hat{K}$} to construct the isoparametric mapping to the physical element $K$. The forcing term of this problem is taken as follows: for $i=Mx$, $f_x = 1$, $f_y = 0$, zero flow rate in the continuity equation and zero traction; for $i=My$, $f_x = 0$, $f_y = 1$, zero flow rate in the continuity equation and zero traction; for $i=C$, $f_x = f_y = 0$, a unit flow rate in the continuity equation and zero traction; for $i=S_{\perp}$, $\bm{f} = \bm{0}$, zero flow rate in the continuity equation and traction of norm 1 parallel to $\mathcal{R}_{S_{\perp}}$; finally, for $i=S_{\parallel}$, $\bm{f} = \bm{0}$, zero flow rate in the continuity equation and traction of norm 1 parallel to $\mathcal{R}_{S_{\parallel}}$. The problems to be solved are of Neumann-type, but the gradient of the velocity solution to these problems is unique. Let us denote by $\nabla\tilde{\bm{u}}_{\mathcal{R}_{i}}$ the gradient of the velocity solution to each of these problems. Since taking a derivative of the velocity will scale like $h_K^{-1}$, these gradients will scale as $\frac{h_K}{\mu}$ for $i = Mx, My$, as a constant for $i = C$ and as $\mu^{-1}$ for $i =S_{\perp}$ and $i=S_{\parallel}$. The way to solve the Stokes problem described is to use as velocity finite element space the one spanned by the interior and edge bubble functions introduced earlier and a pressure space such that the velocity-pressure pair is inf-sup stable. Once this problem is solved, we define
\begin{align}
(\tau)_{\rm visc}^{\mathcal{R}_{i}} = 
\frac{1}{|\hat{K}|}\left(\int_{\hat{K}} 
\nabla\tilde{\bm{u}}_{\mathcal{R}_{i}}:\nabla\tilde{\bm{u}}_{\mathcal{R}_{i}}\right)^{\frac{1}{2}}. \label{eq:tau-visc} 
\end{align}
Since $\hat{K}$ is a fixed parent domain, $(\tau)_{\rm visc}^{\mathcal{R}_{i}}$ will scale like $\nabla\tilde{\bm{u}}_{\mathcal{R}_{i}}$.
Their precise values can be found in \cite{Irisarri2018STOKES} (see Table~B1 in this reference).

Using expressions \eqref{eq:tau-conv} and \eqref{eq:tau-visc} in \eqref{eq:time-scales}, the expression of the SGSs \eqref{eq:GH_Space_Error_Indicator}, and consequently of the local a posteriori error estimator, is now complete.

Once the error for each element is estimated, the global error can be computed as:
\begin{align}
    \eta_{H^1}^{\rm rec} := |\tilde{\bm{u}}_{\rm rec}|_{H^1(\Omega)}=\left(\sum_{K}|\tilde{\bm{u}}_{{\rm rec}}|^{2}_{H^1(K)}\right)^{\frac{1}{2}}.
    \label{eq:Global-H1-Space-error}
\end{align}

\section{Numerical examples}

\noindent In this section, numerical tests are conducted to demonstrate the reliability of the error estimators introduced in this paper. For reference, they are all collected in Table~\ref{tab:error_estimators}, with a brief description of their main properties.

\begin{table}[h] 
\centering
\caption{A posteriori error estimates proposed or evaluated in the paper and their main properties.}
\label{tab:error_estimators}
\renewcommand{\arraystretch}{0.5}
\begin{tabular}{lp{3.0cm}p{9.0cm}}
\hline
\\
\textbf{Estimator} & \textbf{Norm/Target Space} & \textbf{Properties and Key Characteristics} \\ 
\\
\hline
\\

$\eta_{\rm stab}$, Eq.~\eqref{eq:eta-stab}
& Stabilized norm
& 
\begin{minipage}[t]{\linewidth}
\begin{itemize}
    \item Proposed estimator in this paper.
    \item Uses the same SGSs approximation as those employed in the stabilized FE method.
    \item Provides a good and robust a posteriori error estimate only in the stabilized norm.
\end{itemize}
\vspace{1ex}
\end{minipage}
\\
\hline
\\

$\eta_{L^2}$, Eq.~\eqref{eq:eta-L2}
& Velocity $L^2$ norm
& 
\begin{minipage}[t]{\linewidth}

\begin{itemize}
    \item Uses the same SGSs as those employed in the stabilized FE method.
    \item Computes the $L^2$ norm of the SGSs.
    \item There is no theoretical guarantee that it provides a reliable error estimate.
\end{itemize}
\vspace{1ex}
\end{minipage}
\\
\hline
\\

$\eta_{H^1}$, Eq.~\eqref{eq:sgs-eta-h1}
& Velocity $H^1$ seminorm
& 
\begin{minipage}[t]{\linewidth}
\begin{itemize}
    \item Uses the same SGSs as those employed in the stabilized FE method.
    \item Projects the SGSs onto the FE space (smoothing) to compute the $H^1$ seminorm.
    \item There is no theoretical guarantee that it provides a reliable error estimate; moreover, the projection spuriously increases the convergence slope.
\end{itemize}
\vspace{1ex}
\end{minipage}
\\
\hline
\\

$\eta_{H^1}^{\rm rec}$, Eq.~\eqref{eq:Global-H1-Space-error}
& Velocity $H^1$ seminorm
& 
\begin{minipage}[t]{\linewidth}
\begin{itemize}
    \item Reference estimator proposed by Irisarri and Hauke~\cite{Irisarri2021}.
    \item Reconstructs the SGSs by projecting the residuals onto local spaces spanned by
interior and edge bubble functions.
    \item Introduces appropriate scaling factors to compute the $H^1$ velocity seminorm.
\end{itemize}
\vspace{1ex}
\end{minipage}
\\
\hline
\\
\end{tabular}
\end{table}

All quantities are assumed to follow the International System of Units. The first case consists of two manufactured solution tests that help to analyze the behavior of our space error estimator in terms of convergence and effectivity across different flow regimes (Re = 0.01, 1.00, and 200.0). 
Afterward, three other manufactured solutions are employed to extend the study to other norms, such as $L^2$ and $H^1$, and also to compare the proposed estimator with the one proposed in \cite{Irisarri2021}, which is in the $H^1$ seminorm. Additionally, a flow over a backward-facing step for incompressible flow with Reynolds numbers of 1 and 100 is used to compare the space error estimator distribution of the novel approach and that of the reference estimator proposed in Irisarri and Hauke \cite{Irisarri2021}. Finally, to evaluate the effectiveness of the proposed space error estimator within a space refinement strategy relative to a uniform mesh refinement approach, the benchmark problem of incompressible flow over a cylinder at a Reynolds number of 500 is considered.

\subsection{Manufactured solution tests}
\noindent The manufactured solutions are used to study the behavior of the proposed space error estimator. Having a known exact solution is fundamental to quantitatively assess our methodology for error estimation.

In all cases, the prescribed analytical solutions are substituted into the governing equations to derive the associated forcing terms. The computational domain considered for these problems is the unit square $\Omega = (0,1)\times (0,1)\subset \mathbb{R}^2$, with the velocity prescribed on the boundaries of $\Omega$. The domain is partitioned with a structured uniform grid of linear elements of size $h=1/N_p$, where $N_p$ is the number of partitions. The density is set to $\rho = 1$, while various viscosities ($\mu = 100, 1, 0.005$) are employed to represent a wide range of flow regimes. In this study, five manufactured solution pairs are considered $[\bm{u}_{i}(x,y), p_{i}(x,y)]$, $i=1,2,3,4,5$, where $(x,y)$ denotes the Cartesian coordinates, defined as
\begin{align}
\bm{u}_{1}(x,y)
&= 100
[q(x)q'(y),\quad -q'(x)q(y)],
\quad \text{with } q(z)=z^2(1-z)^2
\nonumber\\ 
p_{1}(x,y)
&=100(x-1)^2
\nonumber\\[1.5ex]
\bm{u}_{2}(x,y)
&= 100
[q(x)q'(y),\quad -q'(x)q(y)],
\quad \text{with } q(z)=z^2(1-z)^2\cos(\pi z)
\nonumber\\
p_{2}(x,y)
&=\cos(2x)+\cos(2y)
\nonumber\\[1.5ex]
\bm{u}_{3}(x,y)
&= \bm{w}(x,y,0.1,0.1)\nonumber\\
p_{3}(x,y)
&=r(x,y,0.1,0.1)\nonumber\\[1.5ex]
\bm{u}_{4}(x,y)
&= \bm{w}(x,y,4.0,4.0)\nonumber\\
p_{4}(x,y)
&=r(x,y,4.0,4.0)\nonumber\\[1.5ex]
\bm{u}_{5}(x,y)
&= \bm{w}(x,y,4.2985,0.1)\nonumber\\
p_{5}(x,y)
&=r(x,y,4.2985,0.1)\qquad \text{with }\nonumber\\[1.5ex]
\bm{w}(x,y,R_{1},R_{2})
&=
\left[
\left(1-\cos\left(
2\pi\frac{e^{R_{1}x}-1}{e^{R_{1}}-1}
\right)\right)
\sin\left(
2\pi\frac{e^{R_{2}y}-1}{e^{R_{2}}-1}
\right)
\frac{R_{2}e^{R_{2}y}}
{2\pi(e^{R_{2}}-1)},\quad
\right.
\nonumber\\
&\qquad\left.
-\sin\left(
2\pi\frac{e^{R_{1}x}-1}{e^{R_{1}}-1}
\right)
\left(1-\cos\left(
2\pi\frac{e^{R_{2}y}-1}{e^{R_{2}}-1}
\right)\right)
\frac{R_{1}e^{R_{1}x}}
{2\pi(e^{R_{1}}-1)}
\right]\,,\quad \text{and }
\nonumber\\
r(x,y,R_{1},R_{2})
&=
\frac{
R_{1}R_{2}
\sin\left(
2\pi\frac{e^{R_{1}x}-1}{e^{R_{1}}-1}
\right)
\sin\left(
2\pi\frac{e^{R_{2}y}-1}{e^{R_{2}}-1}
\right)
e^{R_{1}x}e^{R_{2}y}
}{
(e^{R_{1}}-1)(e^{R_{2}}-1)
}
\nonumber
\end{align}

\begin{figure}[htbp]
  \centering
  \begin{minipage}[b]{0.46\textwidth}
    \centering
    \includegraphics[width=\textwidth]{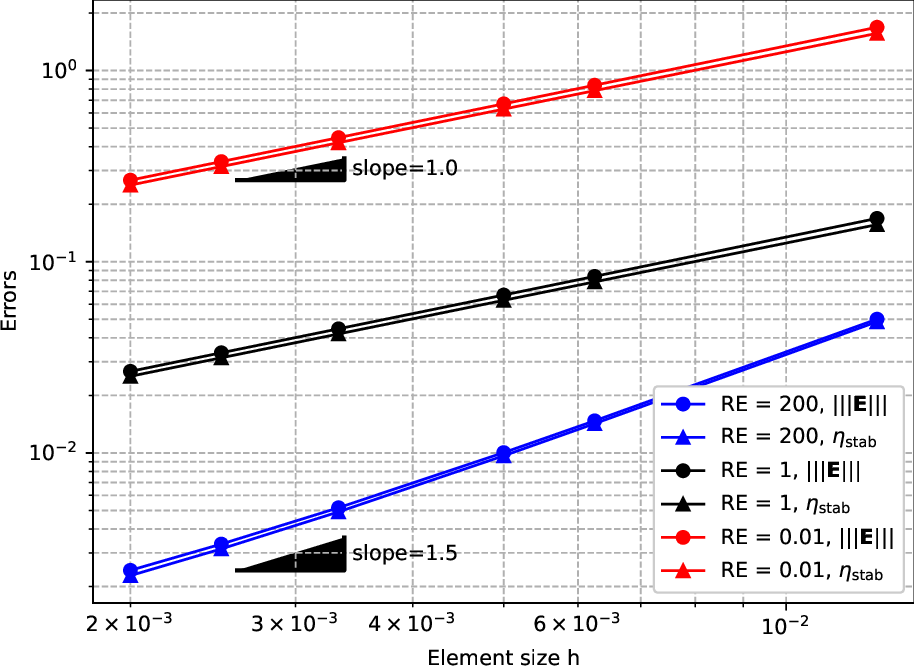}
    \\ \small (a) Errors ($\eta_{\rm stab}$ vs $|||\bm{E}|||$) for $[\bm{u}_{1}, p_{1}]$ .
  \end{minipage}%
  \hfill
  \begin{minipage}[b]{0.46\textwidth}
    \centering
    \includegraphics[width=\textwidth]{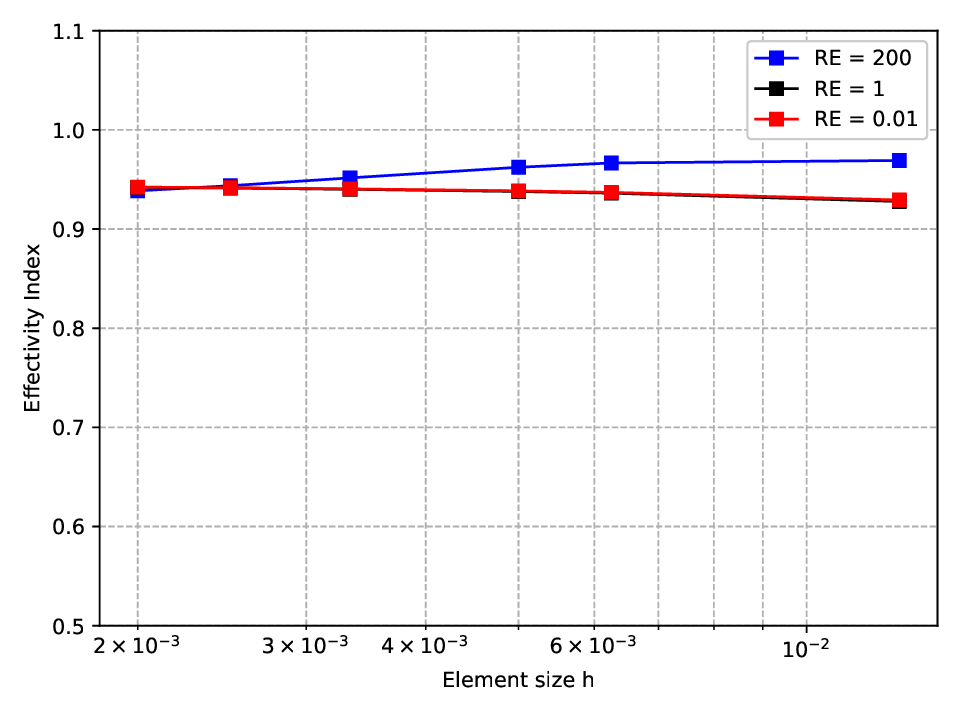}
    \\ \small (b) Effectivity index $(\mathcal{I}_{\rm eff})$ for $[\bm{u}_{1}, p_{1}]$.
  \end{minipage}

  \vspace{1ex} 

  \begin{minipage}[b]{0.46\textwidth}
    \centering
    \includegraphics[width=\textwidth]{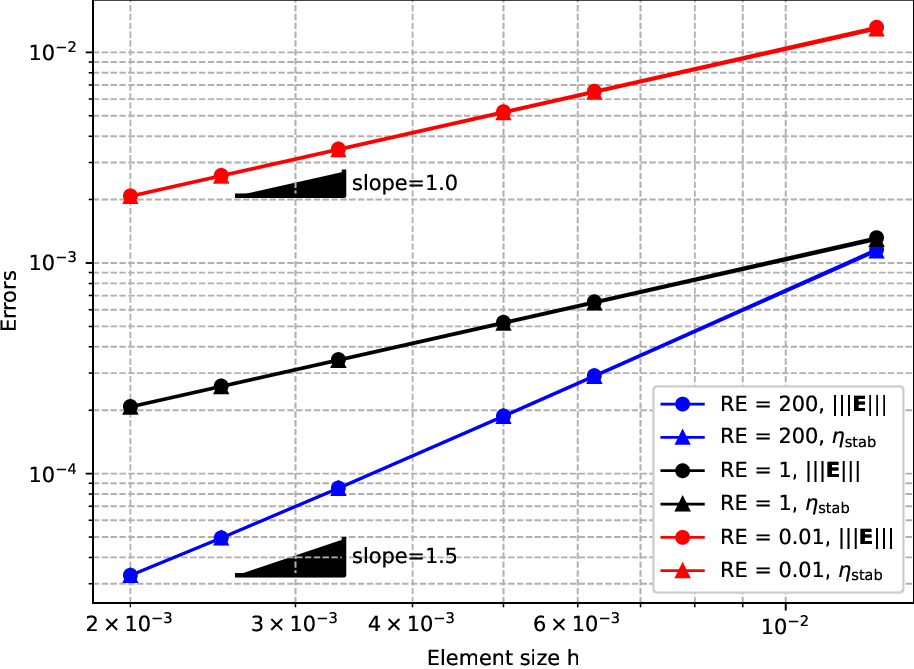}
    \\ \small (c) Errors ($\eta_{\rm stab}$ vs $|||\bm{E}|||$) for $[\bm{u}_{2}, p_{2}]$ .
  \end{minipage}%
  \hfill
  \begin{minipage}[b]{0.46\textwidth}
    \centering
    \includegraphics[width=\textwidth]{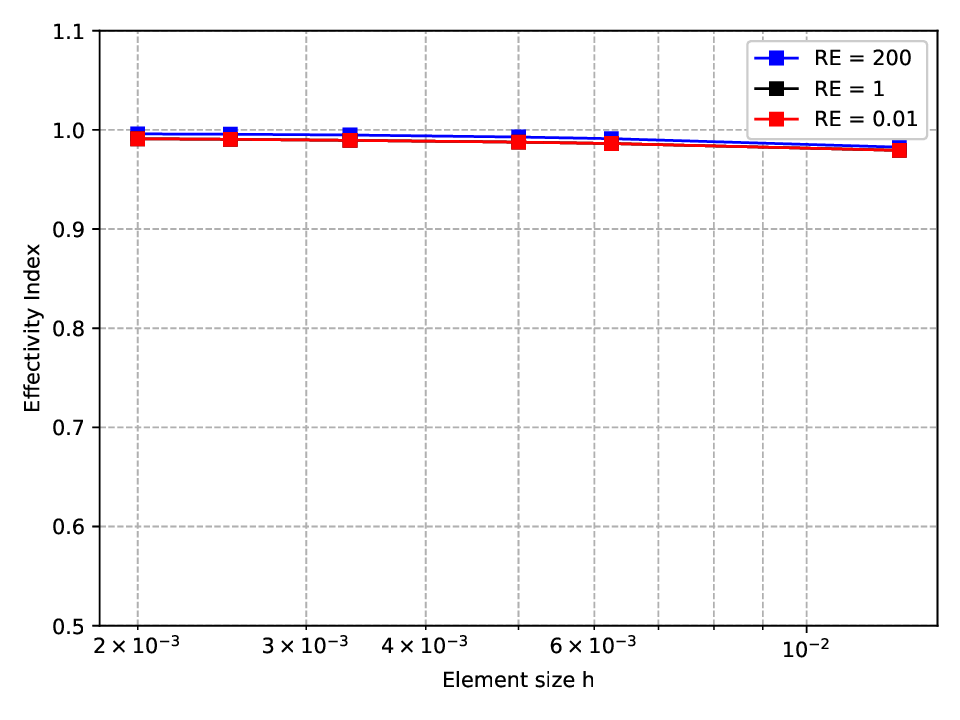}
    \\ \small (d) Effectivity index $(\mathcal{I}_{\rm eff})$ for $[\bm{u}_{2}, p_{2}]$.
  \end{minipage}

  \caption{Space convergence of error and effectivity index results in the stabilized norm ($\eta_{\rm stab}$ vs $|||\bm{E}|||$)  against the element size of linear quadrilateral elements for $\gamma=18.0$, Re = 200, 1, and 0.01 for \textbf{Solution 1} (top) and \textbf{Solution 2} (bottom).}
  \label{fig:part1_STST_Space_Error_Indicator_Cases}
\end{figure}

\begin{figure}[htbp]
  \centering
  \begin{minipage}[b]{0.46\textwidth}
    \centering
    \includegraphics[width=\textwidth]{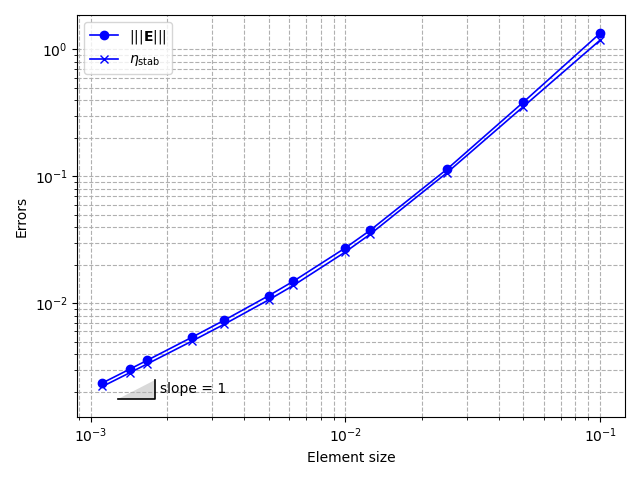}
    \\ 
  \end{minipage}%
  \hfill
  \begin{minipage}[b]{0.46\textwidth}
    \centering
    \includegraphics[width=\textwidth]{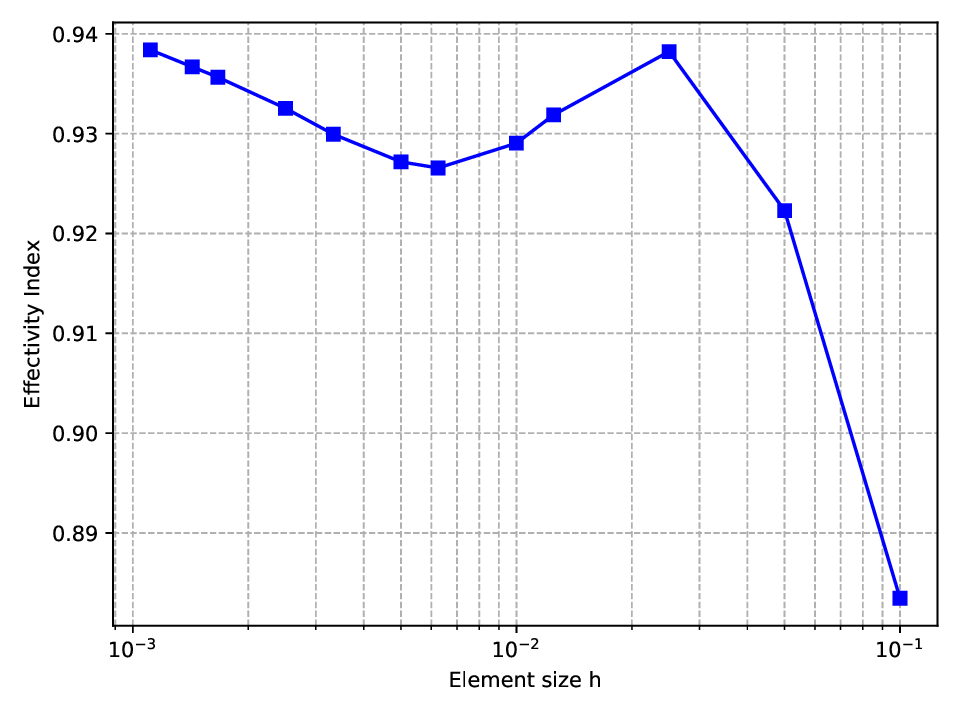}
    \\ 
  \end{minipage}

  \caption{Space convergence of errors ($\eta_{\rm stab}$ vs $|||\bm{E}|||$) on the left and effectivity index on the right $(\mathcal{I}_{\rm eff})$ in the stabilized norm for \textbf{Solution 1}, $[\bm{u}_{1}, p_{1}]$, against the element size of linear quadrilateral elements, $\gamma=18.0$, and Re = 40, transitioning from a convection-dominated to a diffusion-dominated regime.}
  \label{fig:part2_STST_Space_Error_Indicator_Cases}
\end{figure}

Fig.~\ref{fig:part1_STST_Space_Error_Indicator_Cases} presents the results of the mesh convergence study using OSGS for manufactured solutions $[\bm{u}_{i}(x,y), p_{i}(x,y)]$, $i=1,2$. This figure shows the element-size convergence of the SGSs space error estimator alongside the exact error in the stabilized norm for different flow regimes. The Reynolds number, $ \mathrm{Re} =\rho U L / \mu$, is evaluated using a unit characteristic velocity $U=1$, and a unit length $L=1$, and is varied by changing the viscosity $\mu$. These cases are evaluated at Reynolds numbers of 0.01, 1.0, and 200.0. As seen from the convergence plots, the proposed SGS-based error estimator matches the behavior of the stabilized norm of the exact error for all three flow regimes. In terms of convergence rate, for the convection-dominated cases ($\mathrm{Re}=200$), the convergence slope of the error estimator is 1.5, compared to viscous-dominated cases where the slope changes to 1. This is consistent with theoretical expectations by looking at the dominant terms of the estimator. The slope transition from convective to viscous-dominated is visible in Fig.~\ref{fig:part2_STST_Space_Error_Indicator_Cases}, where the overall Reynolds number is constant, but as the mesh is continuously refined, the local Reynolds number decreases, and the viscous term starts to dominate. When comparing these results to those from \cite{codina2025posteriorierrorestimatesfinite}, it is observed that, to achieve an effectivity index close to 1, especially in viscous-dominated cases, it is necessary to use $\gamma \approx 18.0$. The reason is that, in the space error estimator, an approximation (\ref{eq:equiv}) is used instead of the exact term for the velocity SGSs on the element boundaries. This approximation requires tuning the constant $\gamma$. As shown in these two figures, across various tested flow regimes and flow patterns, the effectivity index exceeded 0.935, indicating strong alignment between the error estimator and the exact error.

\begin{figure}[htbp]
  \centering
  \begin{minipage}[b]{0.46\textwidth}
    \centering
    \includegraphics[width=\textwidth]{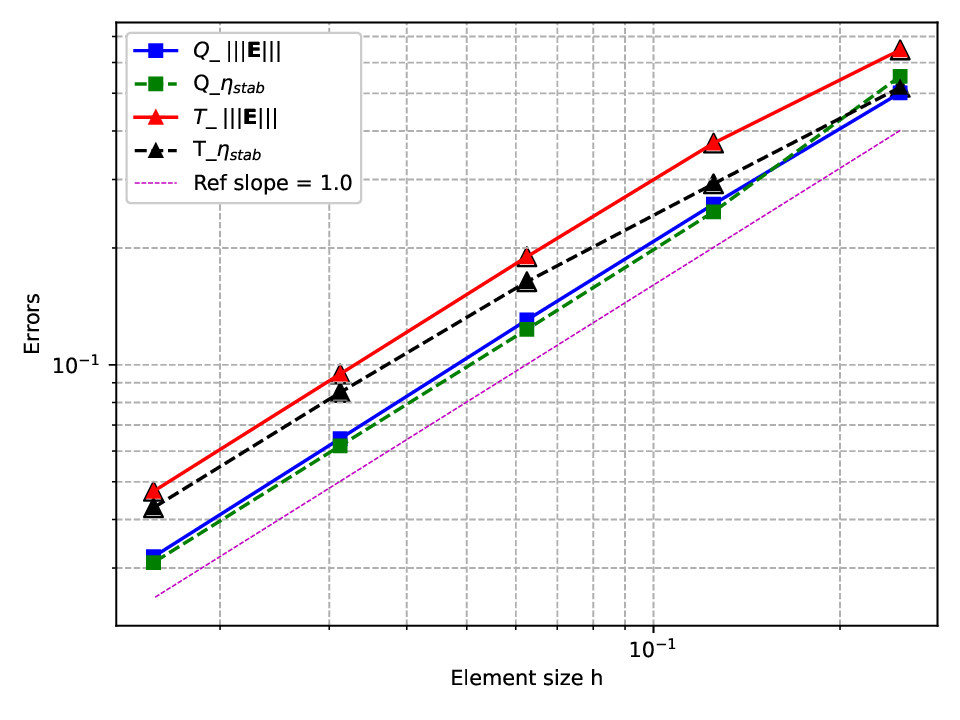}
    \\ \small (a) Stabilized norm ($\eta_{\rm stab}$ vs $|||\bm{E}|||$).
  \end{minipage}%
  \hfill
  \begin{minipage}[b]{0.46\textwidth}
    \centering
    \includegraphics[width=\textwidth]{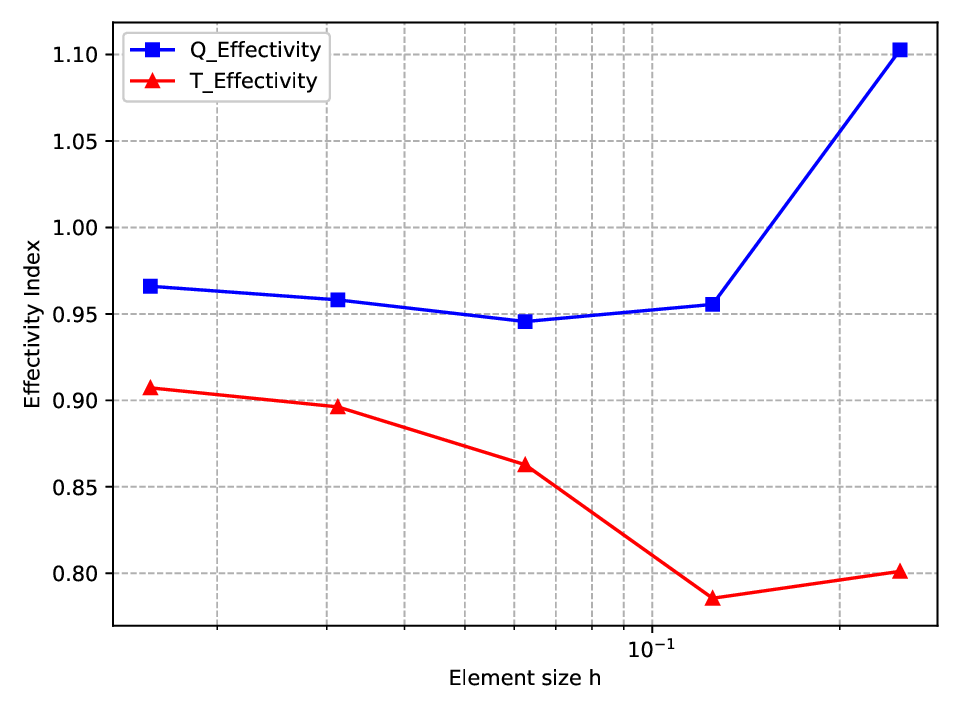}
    \\ \small (b) Effectivity index in stabilized norm.
  \end{minipage}

  \vspace{1ex} 

  \begin{minipage}[b]{0.46\textwidth}
    \centering
    \includegraphics[width=\textwidth]{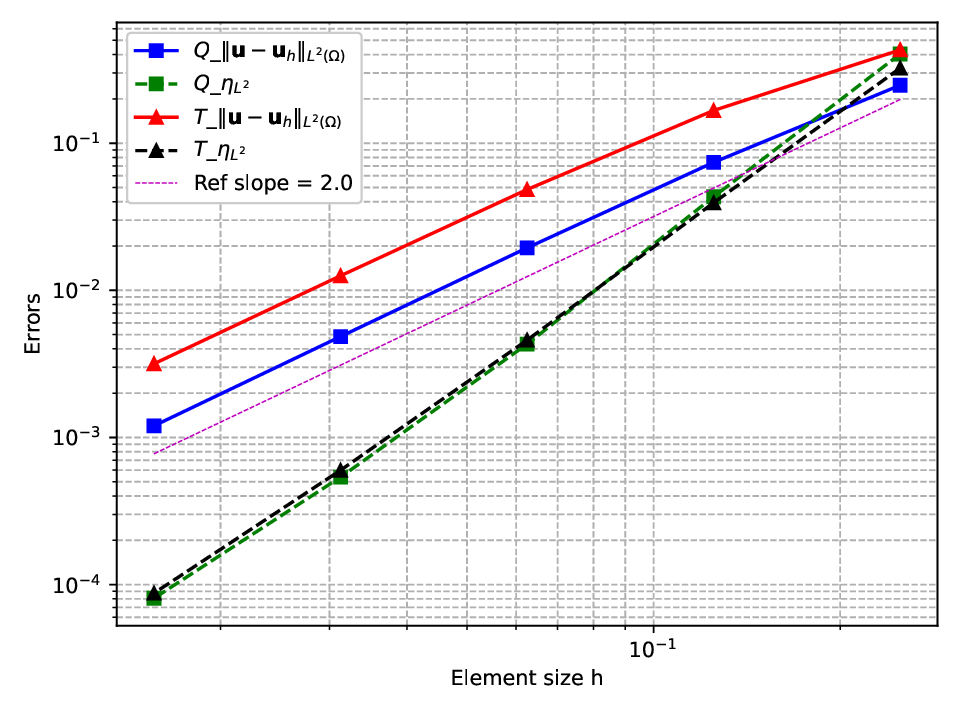}
    \\ \small (c)  $L^2$ norm ($\eta_{L^2}$ vs $E_{L^2}$).
  \end{minipage}%
  \hfill
  \begin{minipage}[b]{0.46\textwidth}
    \centering
    \includegraphics[width=\textwidth]{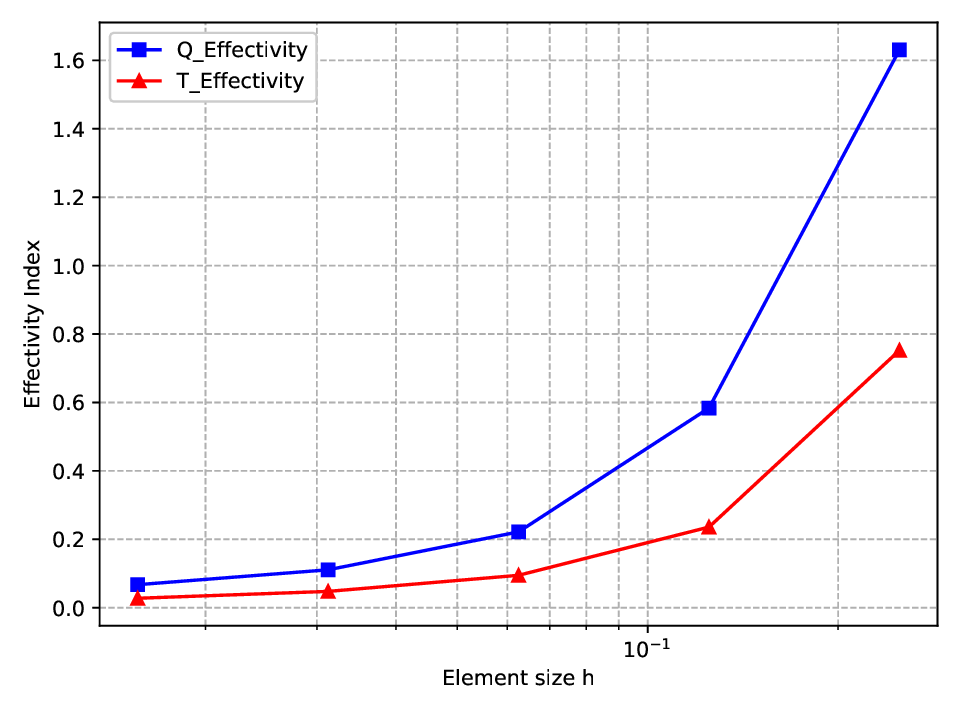}
    \\ \small (d) Effectivity index in $L^2$ norm.
  \end{minipage}

  \begin{minipage}[b]{0.46\textwidth}
    \centering
    \includegraphics[width=\textwidth]{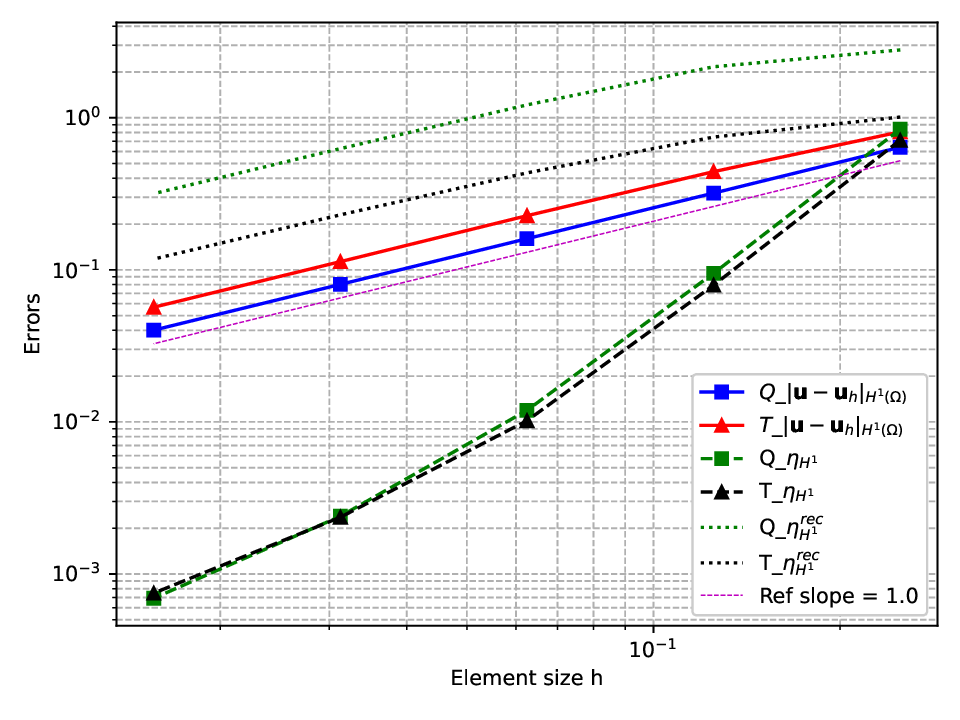}
    \\ \small (e) $H^1$ seminorm ($\eta_{H^1}$, $\eta_{H^1}^{\rm rec}$ vs $E_{H^1}$).
  \end{minipage}
  \hfill
  \begin{minipage}[b]{0.46\textwidth}
    \centering
    \includegraphics[width=\textwidth]{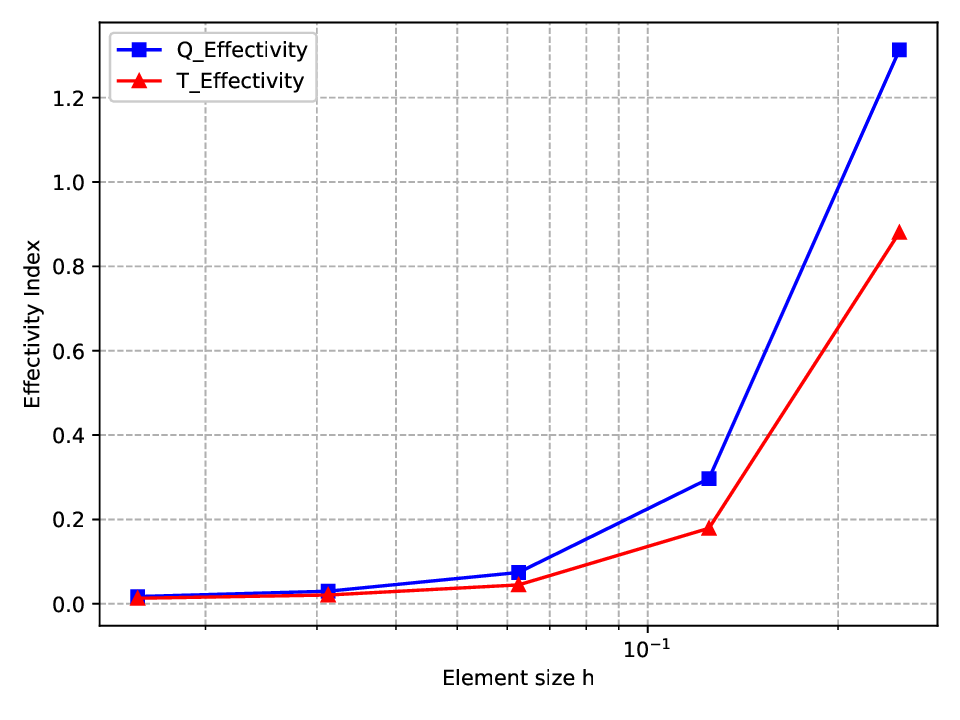}
    \\ \small (f) Effectivity index in $H^1$ seminorm.
  \end{minipage}
  
  \caption{Space convergence (on the left) and effectivity index results (on the right) against the element size of quadrilateral (Q) and triangular (T) linear elements for \textbf{solution 3}  $[\bm{u}_3, p_3]$, $\gamma=18.0$, $\mu  = 1$, and \textbf{OSGS} stabilization technique.}
  \label{fig:OSGS_Test_case_1}
\end{figure}

\begin{figure}[htbp]
  \centering
  \begin{minipage}[b]{0.46\textwidth}
    \centering
    \includegraphics[width=\textwidth]{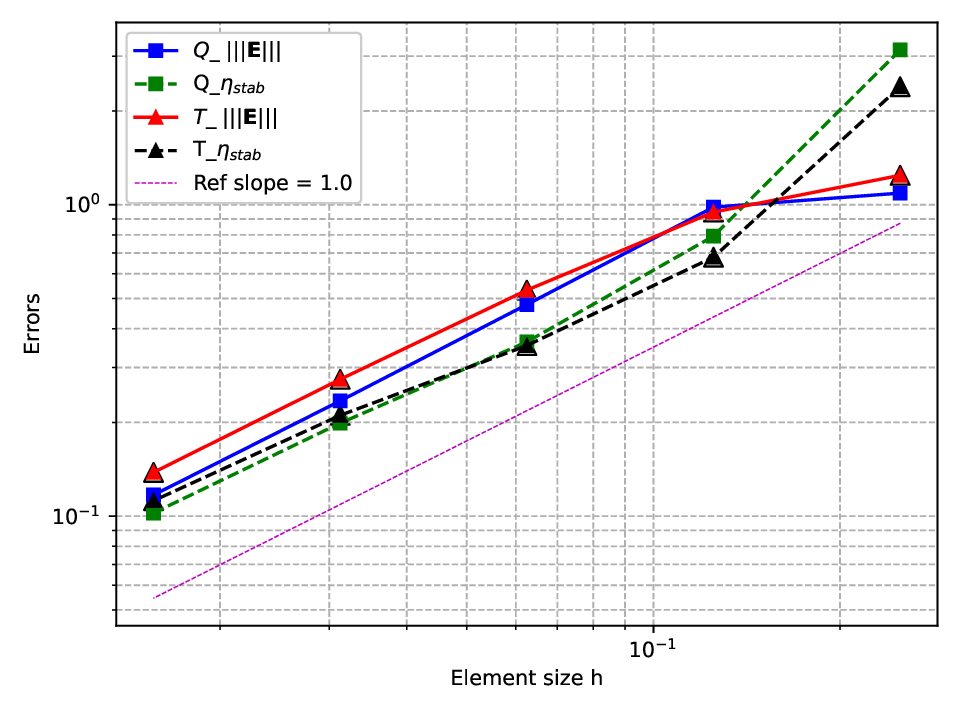}
    \\ \small (a) Stabilized norm ($\eta_{\rm stab}$ vs $|||\bm{E}|||$).
  \end{minipage}%
  \hfill
  \begin{minipage}[b]{0.46\textwidth}
    \centering
    \includegraphics[width=\textwidth]{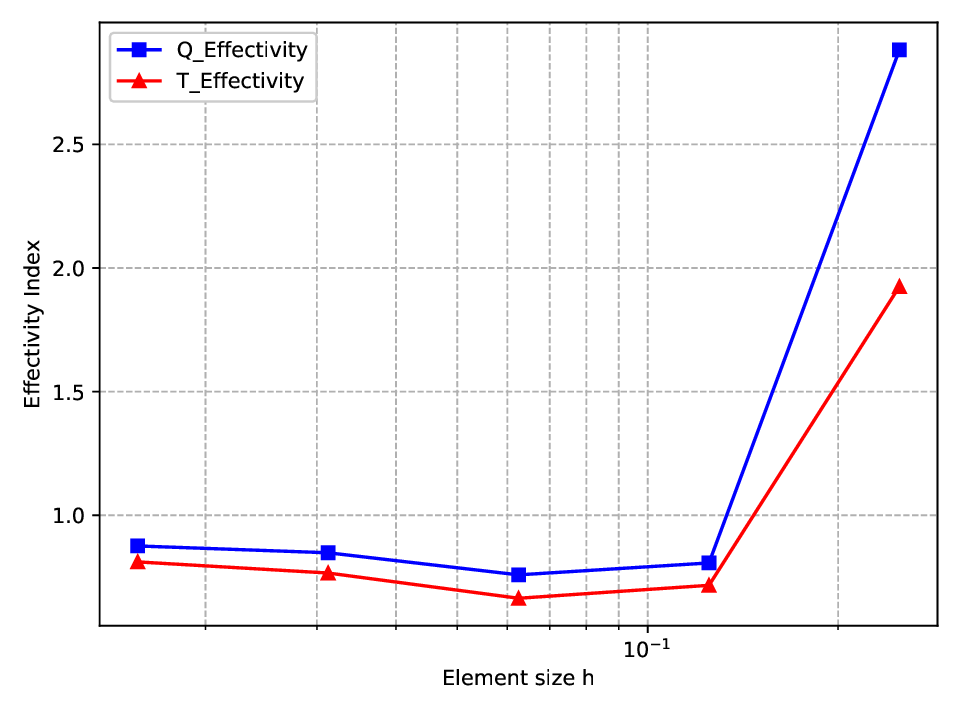}
    \\ \small (b) Effectivity index in stabilized norm.
  \end{minipage}

  \vspace{1ex} 

  \begin{minipage}[b]{0.46\textwidth}
    \centering
    \includegraphics[width=\textwidth]{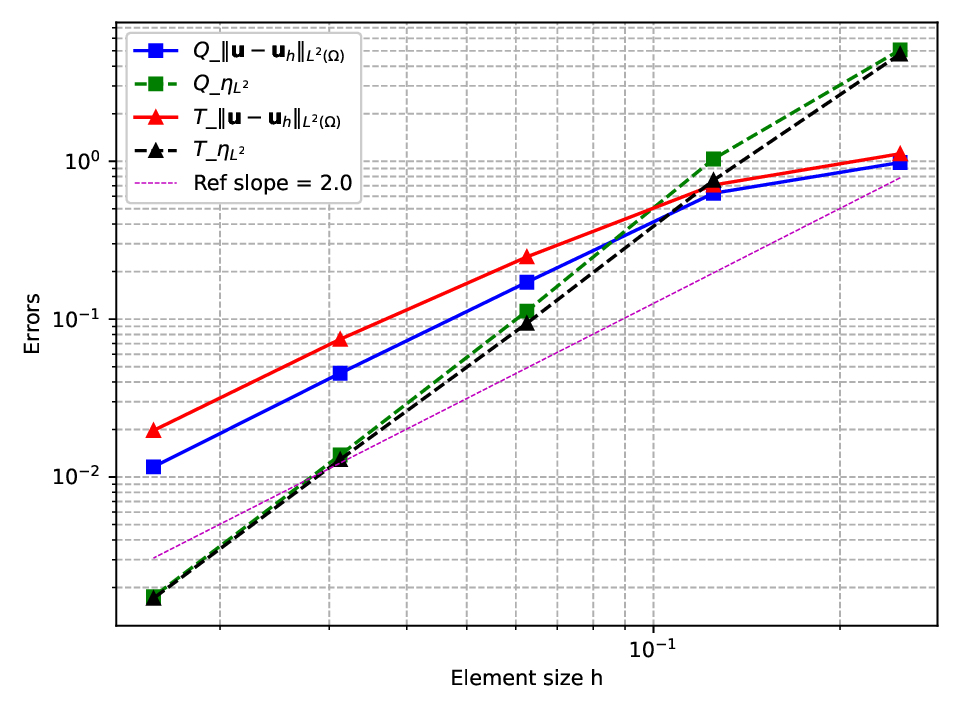}
    \\ \small (c)  $L^2$ norm ($\eta_{L^2}$ vs $E_{L^2}$).
  \end{minipage}%
  \hfill
  \begin{minipage}[b]{0.46\textwidth}
    \centering
    \includegraphics[width=\textwidth]{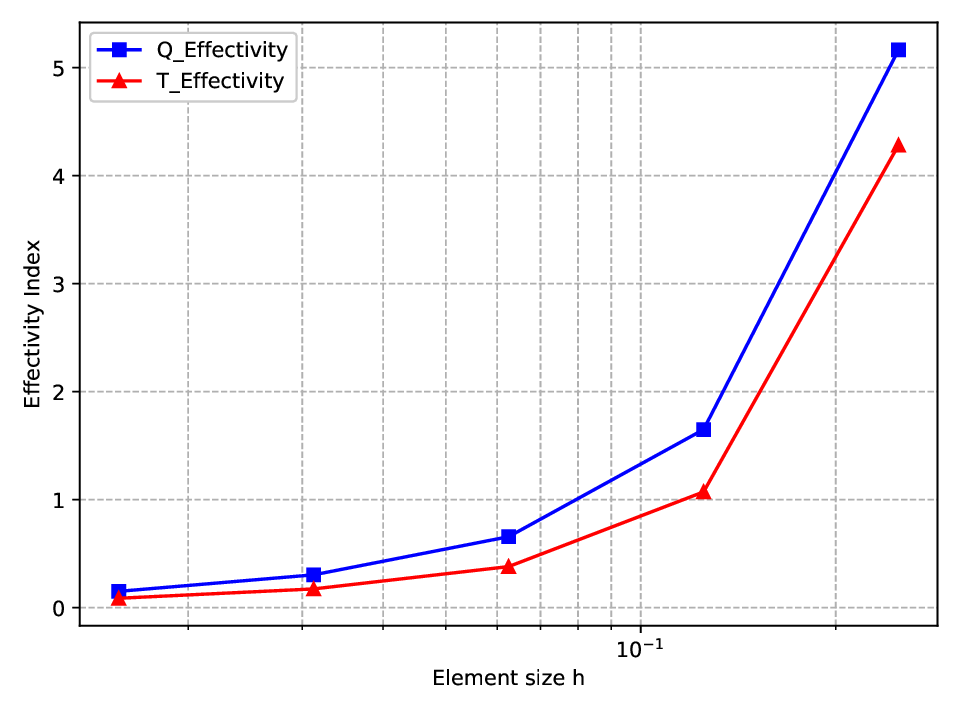}
    \\ \small (d) Effectivity index in $L^2$ norm.
  \end{minipage}

  \begin{minipage}[b]{0.46\textwidth}
    \centering
    \includegraphics[width=\textwidth]{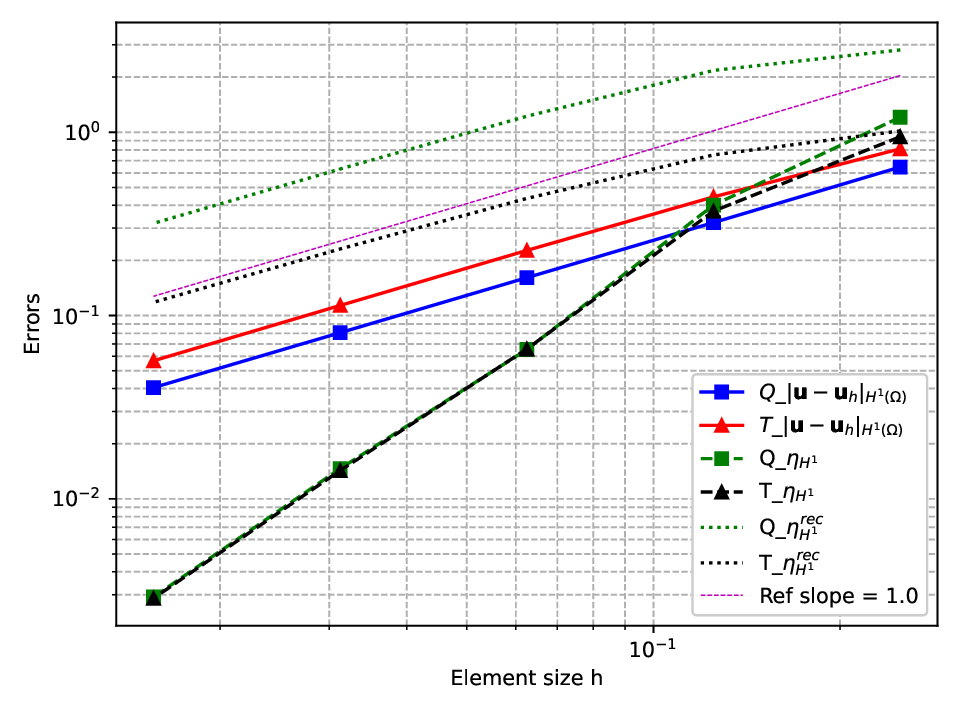}
    \\ \small (e) $H^1$ seminorm ($\eta_{H^1}$, $\eta_{H^1}^{\rm rec}$ vs $E_{H^1}$).
  \end{minipage}%
  \hfill
  \begin{minipage}[b]{0.46\textwidth}
    \centering
    \includegraphics[width=\textwidth]{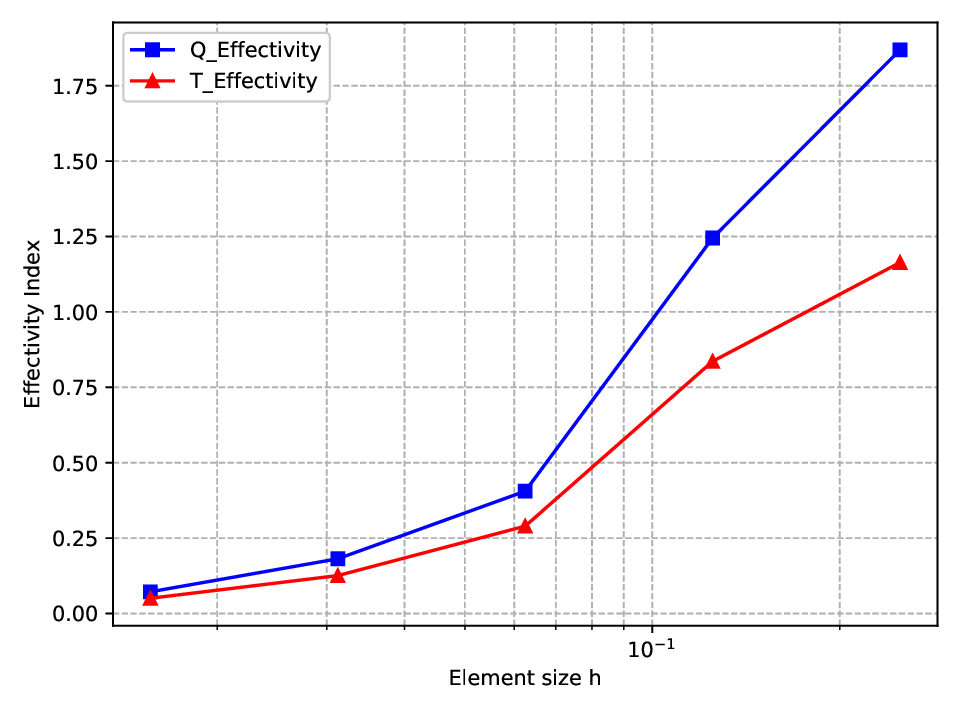}
    \\ \small (f) Effectivity index in $H^1$ seminorm.
  \end{minipage}
  
  \caption{Space convergence (on the left) and effectivity index results (on the right) against the element size of quadrilateral (Q) and triangular (T) linear elements for \textbf{solution 4}  $[\bm{u}_4, p_4]$, $\gamma=18.0$, $\mu  = 1$, and \textbf{OSGS} stabilization technique.}
  \label{fig:OSGS_Test_case_2}
\end{figure}

\begin{figure}[htbp]
  \centering
  \begin{minipage}[b]{0.46\textwidth}
    \centering
    \includegraphics[width=\textwidth]{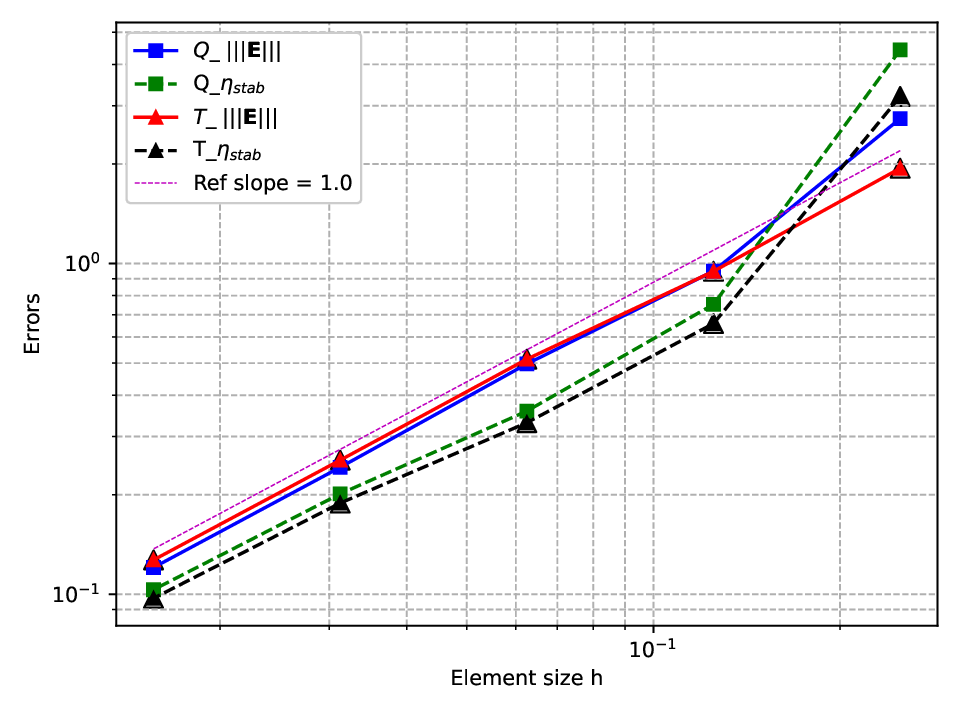}
    \\ \small (a) Stabilized norm ($\eta_{\rm stab}$ vs $|||\bm{E}|||$).
  \end{minipage}%
  \hfill
  \begin{minipage}[b]{0.46\textwidth}
    \centering
    \includegraphics[width=\textwidth]{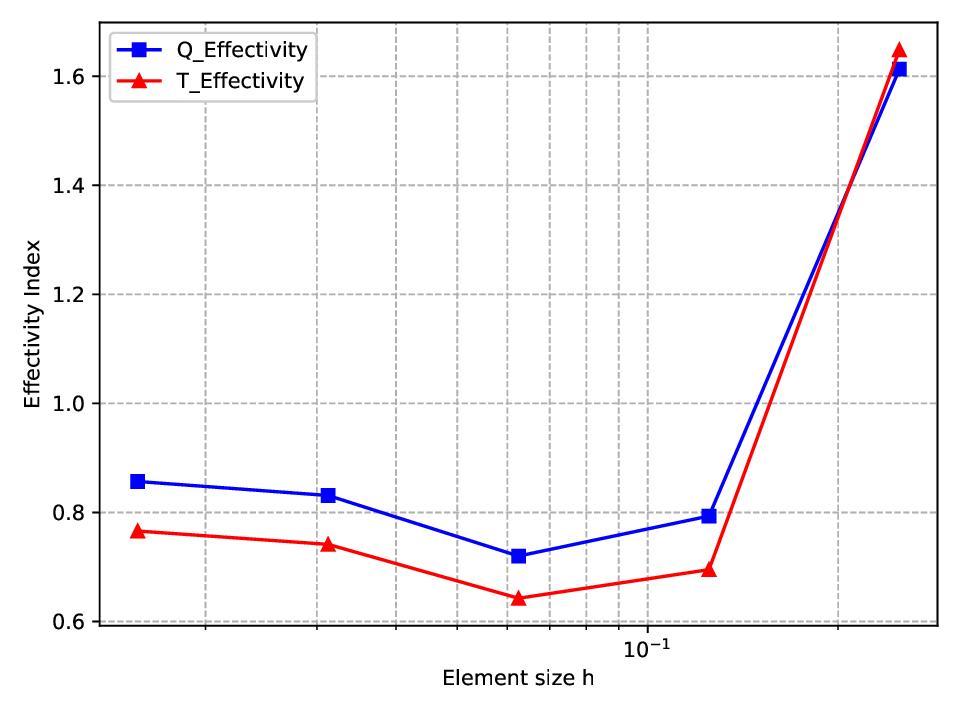}
    \\ \small (b) Effectivity index in stabilized norm.
  \end{minipage}

  \vspace{1ex} 

  \begin{minipage}[b]{0.46\textwidth}
    \centering
    \includegraphics[width=\textwidth]{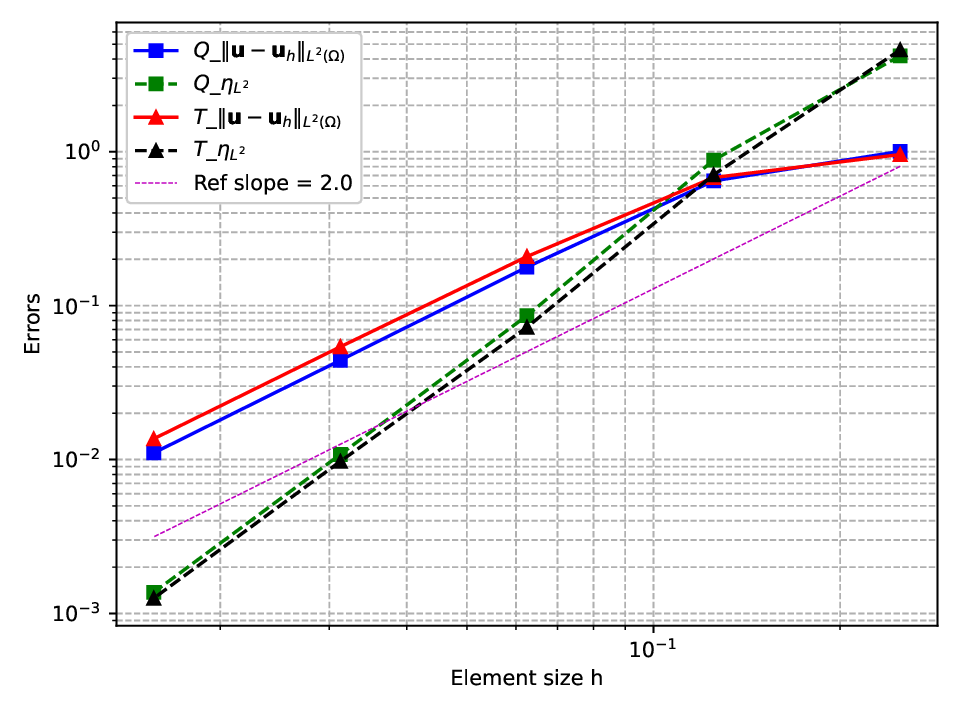}
    \\ \small (c)  $L^2$ norm ($\eta_{L^2}$ vs $E_{L^2}$).
  \end{minipage}%
  \hfill
  \begin{minipage}[b]{0.46\textwidth}
    \centering
    \includegraphics[width=\textwidth]{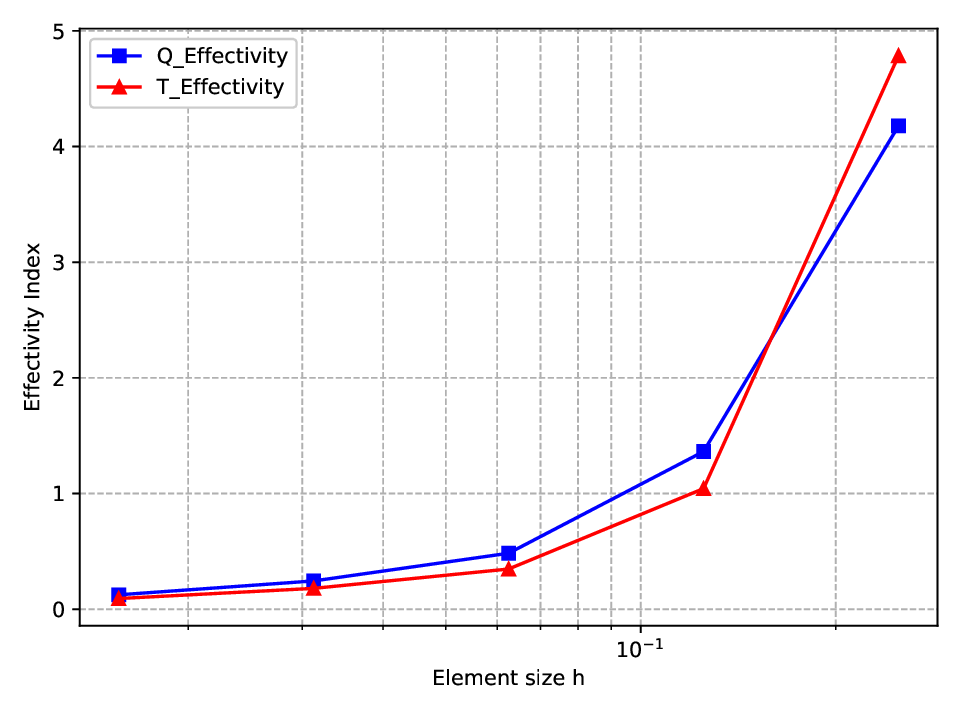}
    \\ \small (d) Effectivity index in $L^2$ norm.
  \end{minipage}

  \begin{minipage}[b]{0.46\textwidth}
    \centering
    \includegraphics[width=\textwidth]{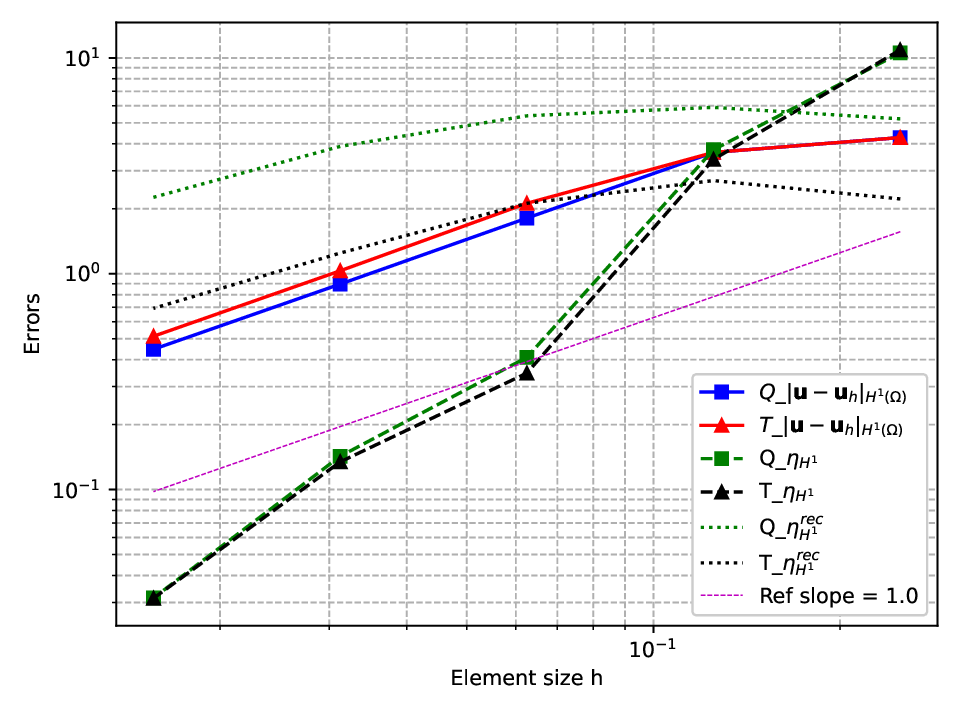}
    \\ \small (e) $H^1$ seminorm ($\eta_{H^1}$, $\eta_{H^1}^{\rm rec}$ vs $E_{H^1}$).
  \end{minipage}%
  \hfill
  \begin{minipage}[b]{0.46\textwidth}
    \centering
    \includegraphics[width=\textwidth]{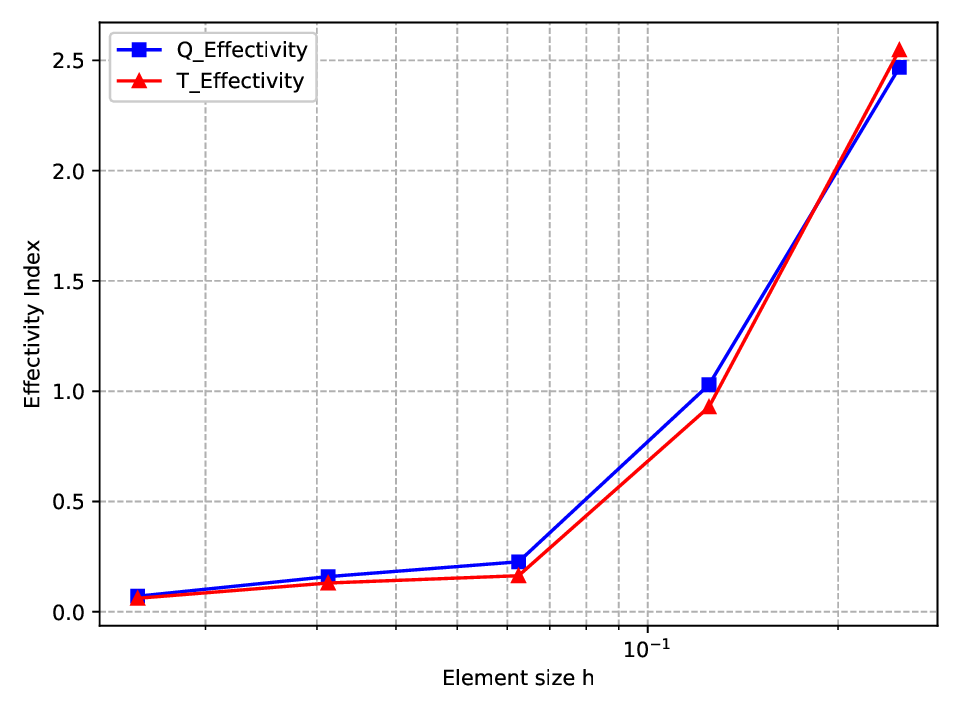}
    \\ \small (f) Effectivity index in $H^1$ seminorm.
  \end{minipage}
  
  \caption{Space convergence (on the left) and effectivity index results (on the right) against the element size of quadrilateral (Q) and triangular (T) linear elements for \textbf{solution 5}  $[\bm{u}_5, p_5]$, $\gamma=18.0$, $\mu  = 1$, and \textbf{OSGS} stabilization technique.}
  \label{fig:OSGS_Test_case_3}
\end{figure}

\begin{figure}[htbp]
  \centering
  \begin{minipage}[b]{0.46\textwidth}
    \centering
    \includegraphics[width=\textwidth]{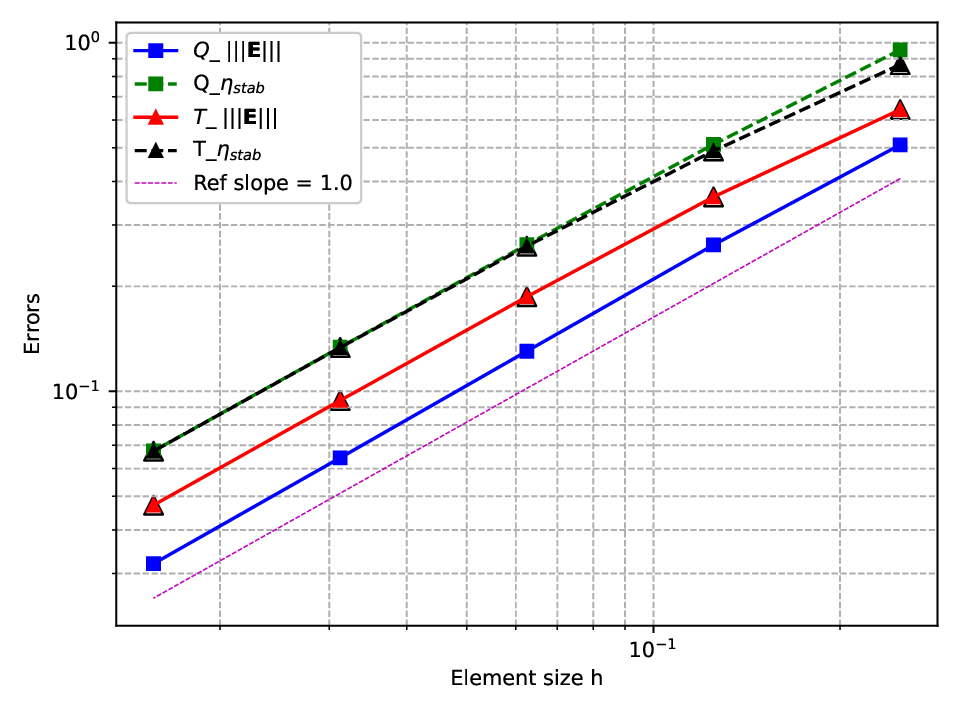}
    \\ \small (a) Stabilized norm ($\eta_{\rm stab}$ vs $|||\bm{E}|||$).
  \end{minipage}%
  \hfill
  \begin{minipage}[b]{0.46\textwidth}
    \centering
    \includegraphics[width=\textwidth]{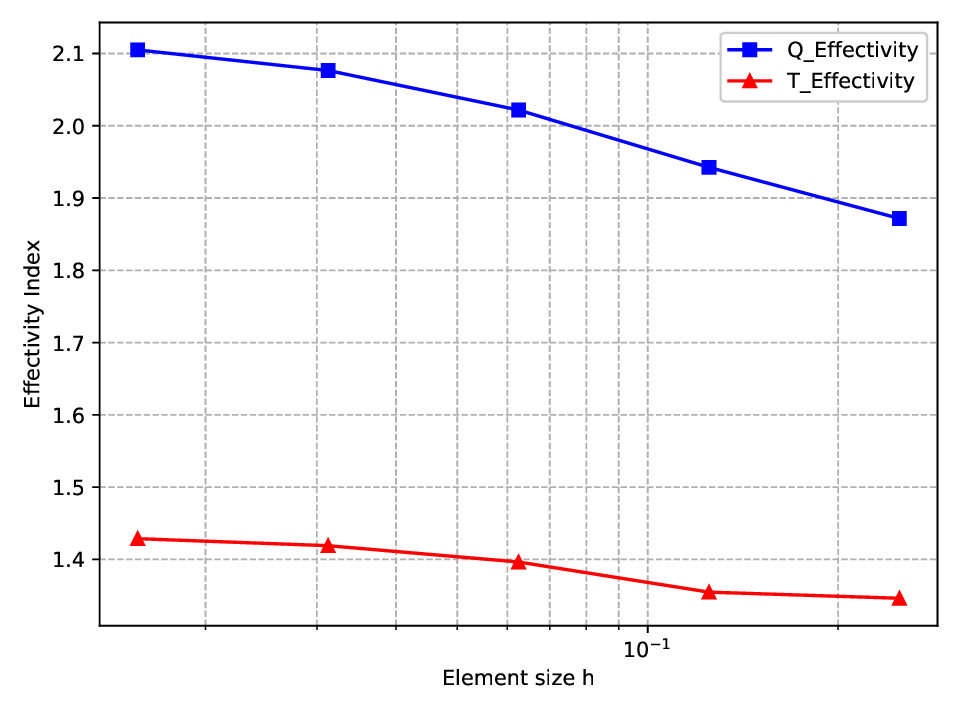}
    \\ \small (b) Effectivity index in stabilized norm.
  \end{minipage}

  \vspace{1ex} 

  \begin{minipage}[b]{0.46\textwidth}
    \centering
    \includegraphics[width=\textwidth]{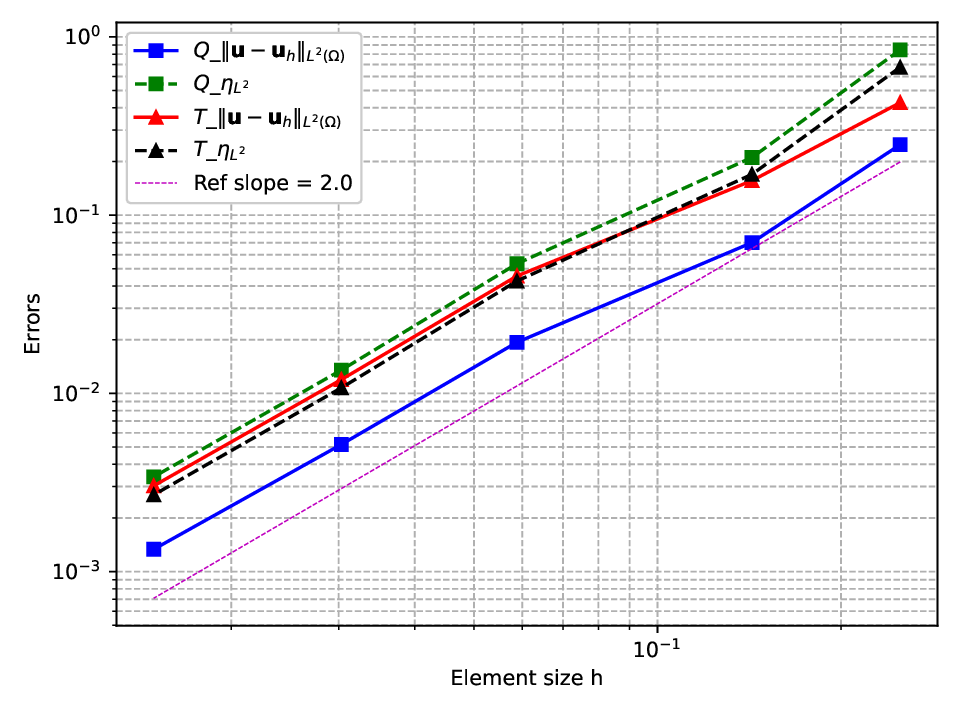}
    \\ \small (c)  $L^2$ norm ($\eta_{L^2}$ vs $E_{L^2}$).
  \end{minipage}%
  \hfill
  \begin{minipage}[b]{0.46\textwidth}
    \centering
    \includegraphics[width=\textwidth]{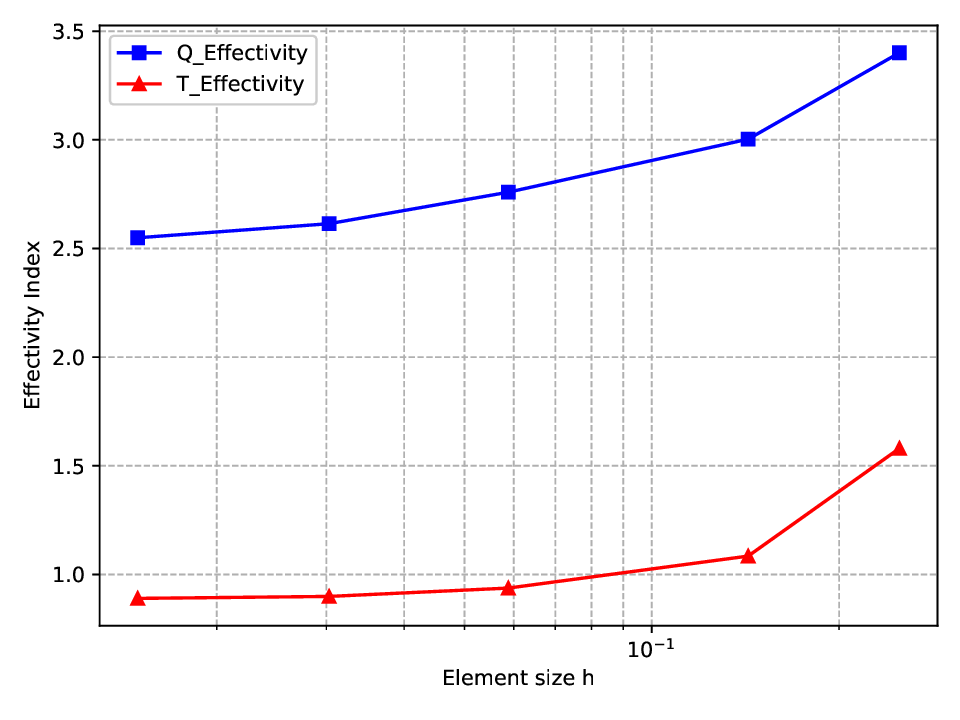}
    \\ \small (d) Effectivity index in $L^2$ norm.
  \end{minipage}

  \begin{minipage}[b]{0.46\textwidth}
    \centering
    \includegraphics[width=\textwidth]{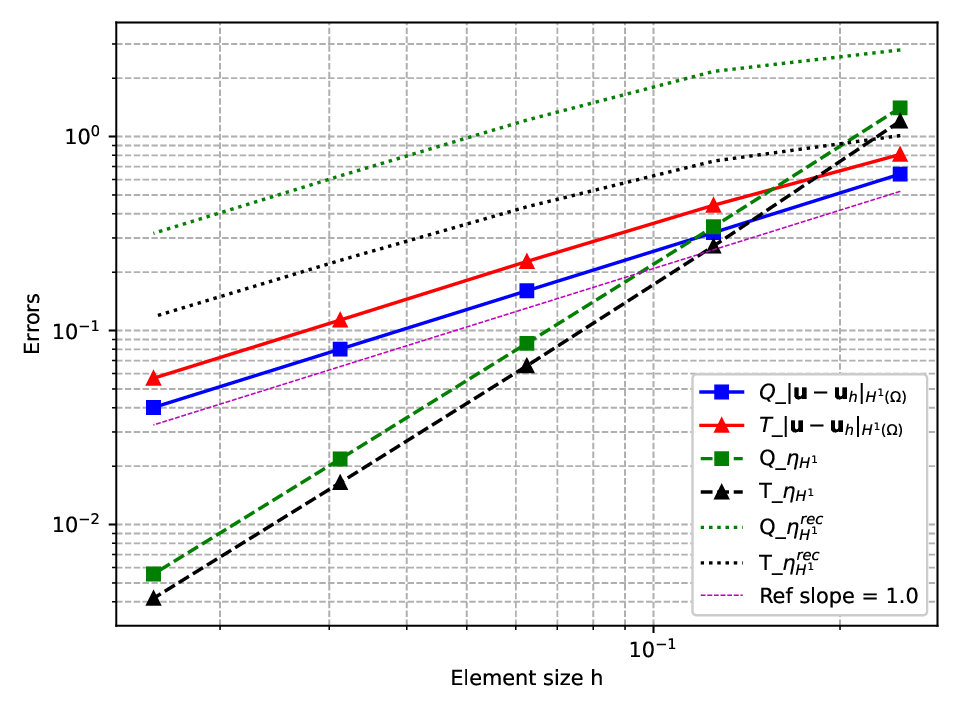}
    \\ \small (e) $H^1$ seminorm ($\eta_{H^1}$, $\eta_{H^1}^{\rm rec}$ vs $E_{H^1}$).
  \end{minipage}%
  \hfill
  \begin{minipage}[b]{0.46\textwidth}
    \centering
    \includegraphics[width=\textwidth]{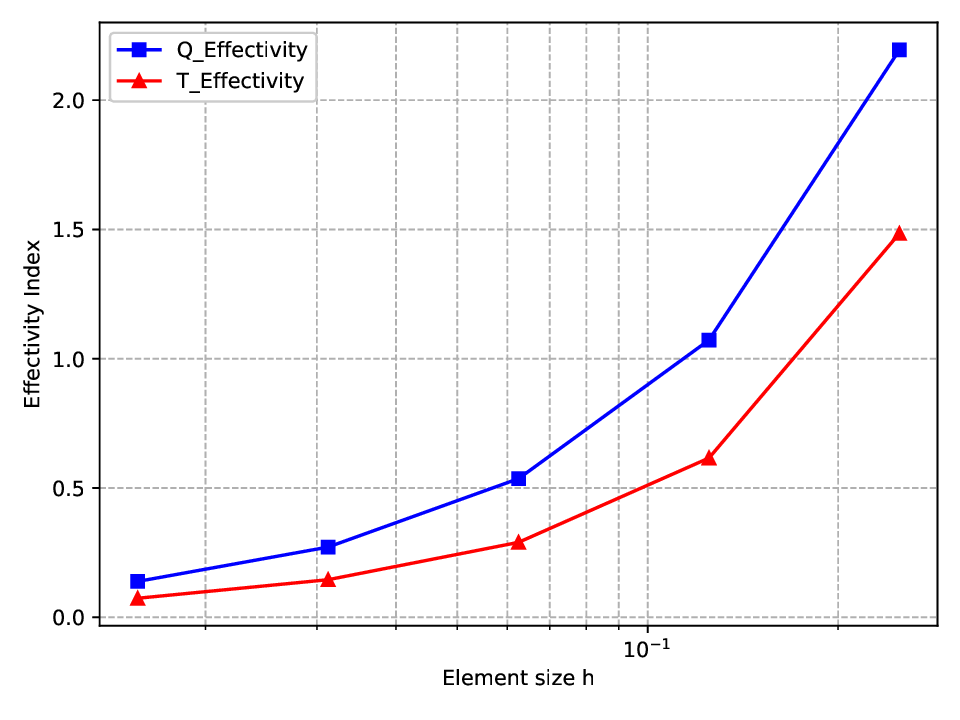}
    \\ \small (f) Effectivity index in $H^1$ seminorm.
  \end{minipage}
  
  \caption{Space convergence (on the left) and effectivity index results (on the right) against the element size of quadrilateral (Q) and triangular (T) linear elements for \textbf{solution 3}  $[\bm{u}_3, p_3]$, $\gamma=18.0$, $\mu  = 1$, and \textbf{ASGS} stabilization technique.}
  \label{fig:ASGS_Test_case_1}
\end{figure}

\begin{figure}[htbp]
  \centering
  \begin{minipage}[b]{0.46\textwidth}
    \centering
    \includegraphics[width=\textwidth]{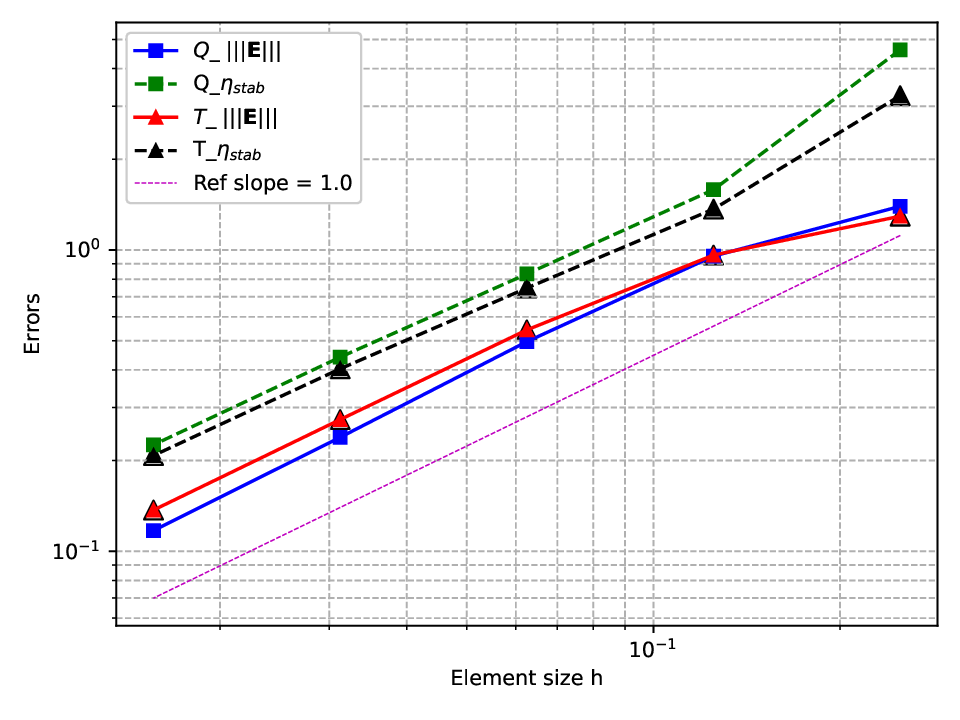}
    \\ \small (a) Stabilized norm ($\eta_{\rm stab}$ vs $|||\bm{E}|||$).
  \end{minipage}%
  \hfill
  \begin{minipage}[b]{0.46\textwidth}
    \centering
    \includegraphics[width=\textwidth]{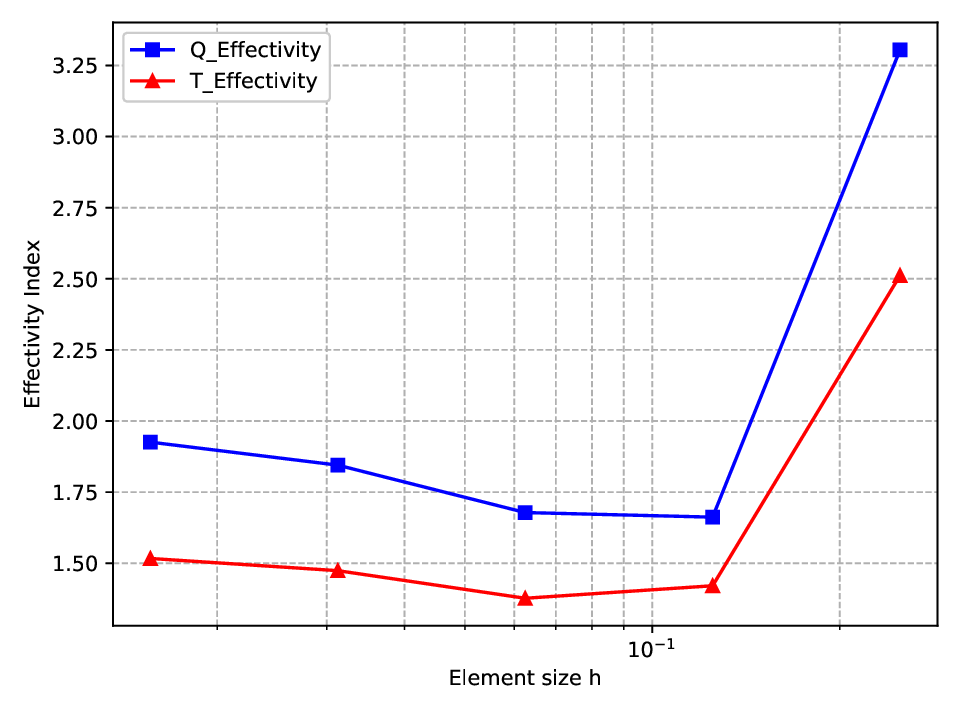}
    \\ \small (b) Effectivity index in stabilized norm.
  \end{minipage}

  \vspace{1ex} 

  \begin{minipage}[b]{0.46\textwidth}
    \centering
    \includegraphics[width=\textwidth]{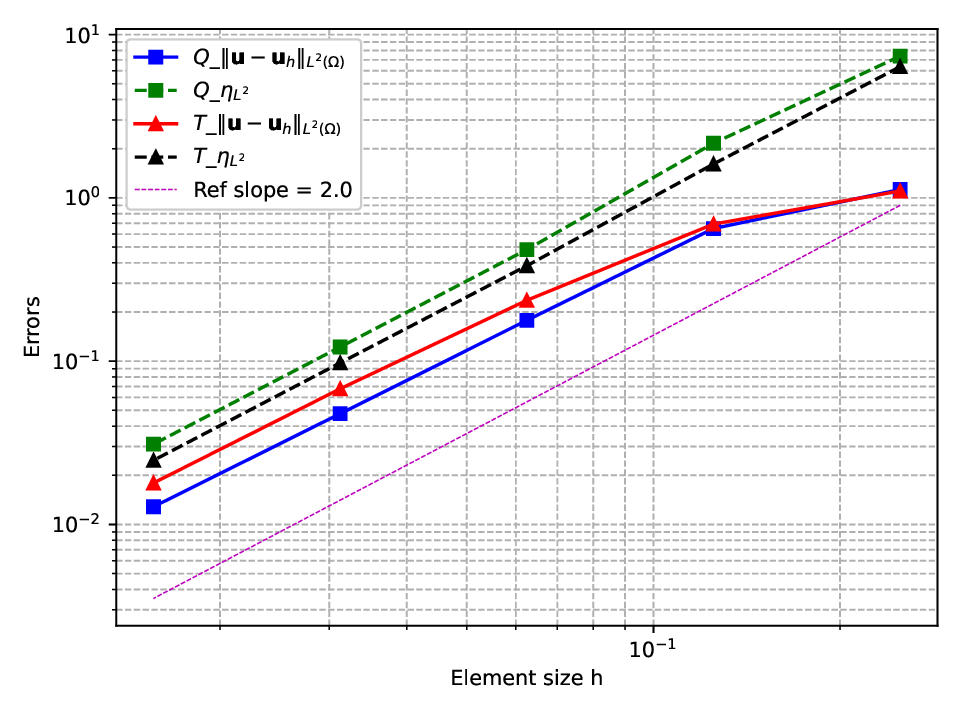}
    \\ \small (c)  $L^2$ norm ($\eta_{L^2}$ vs $E_{L^2}$).
  \end{minipage}%
  \hfill
  \begin{minipage}[b]{0.46\textwidth}
    \centering
    \includegraphics[width=\textwidth]{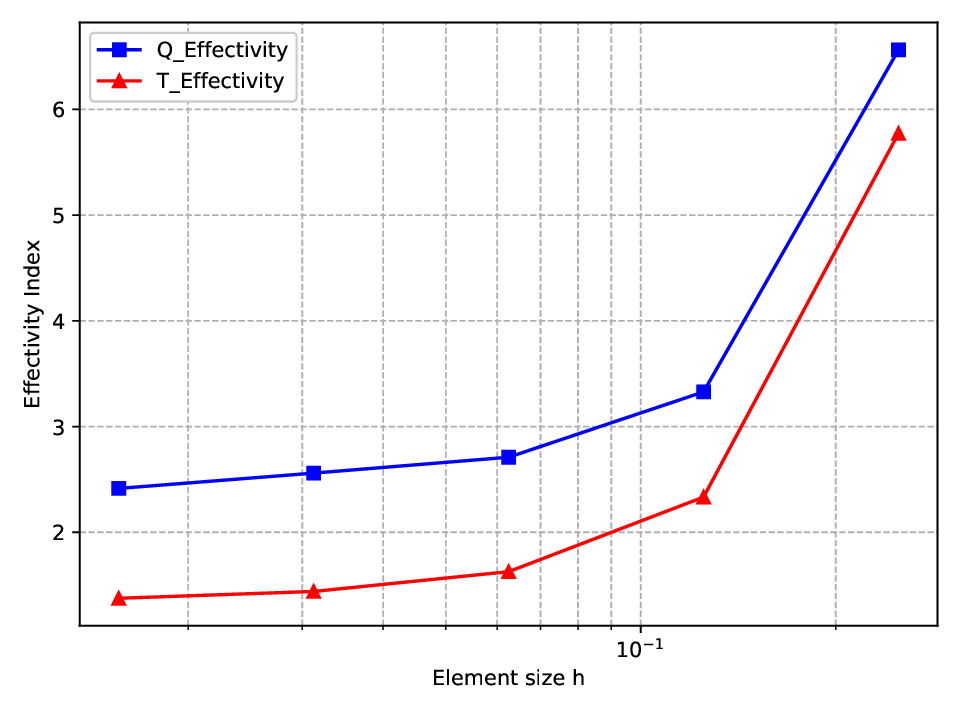}
    \\ \small (d) Effectivity index in $L^2$ norm.
  \end{minipage}

  \begin{minipage}[b]{0.46\textwidth}
    \centering
    \includegraphics[width=\textwidth]{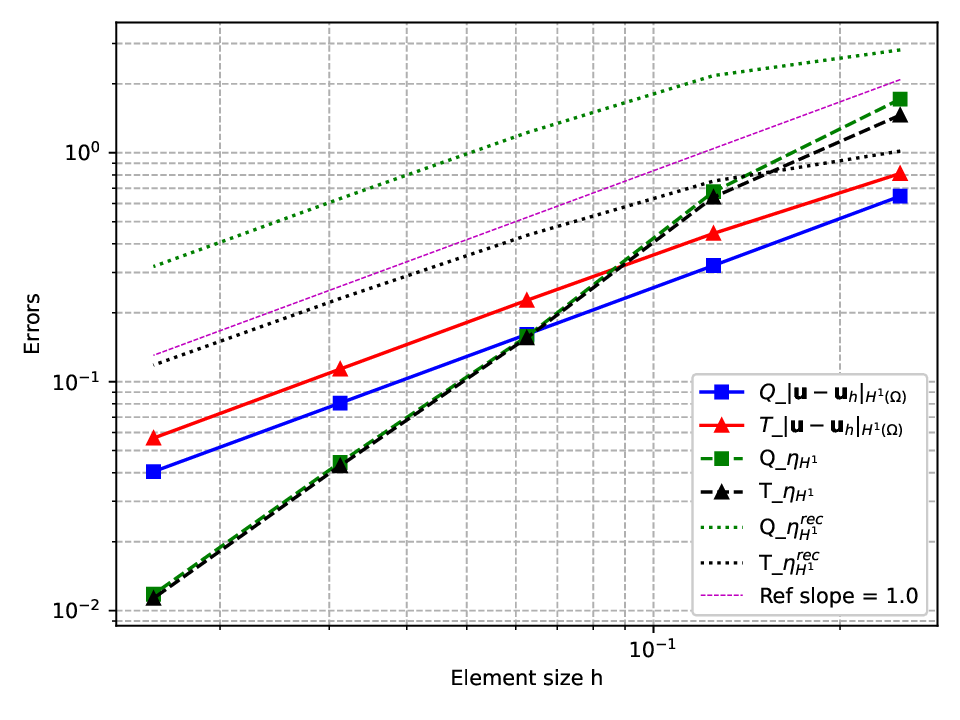}
    \\ \small (e) $H^1$ seminorm ($\eta_{H^1}$, $\eta_{H^1}^{\rm rec}$ vs $E_{H^1}$).
  \end{minipage}%
  \hfill
  \begin{minipage}[b]{0.46\textwidth}
    \centering
    \includegraphics[width=\textwidth]{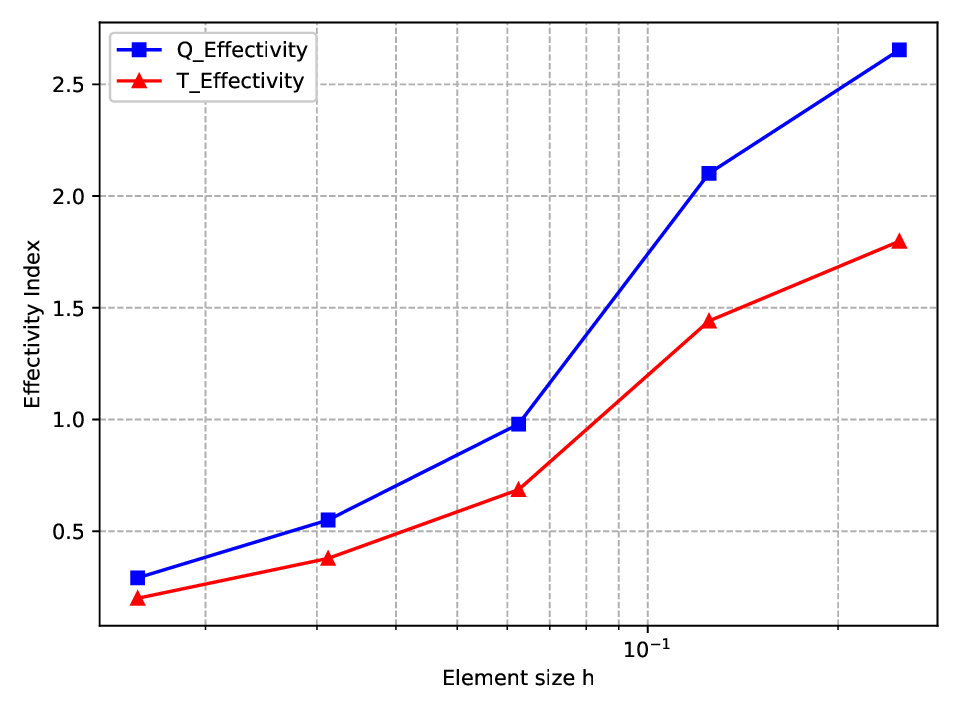}
    \\ \small (f) Effectivity index in $H^1$ seminorm.
  \end{minipage}
  
  \caption{Space convergence (on the left) and effectivity index results (on the right) against the element size of quadrilateral (Q) and triangular (T) linear elements for \textbf{solution 4}  $[\bm{u}_4, p_4]$, $\gamma=18.0$, $\mu  = 1$, and \textbf{ASGS} stabilization technique.}
  \label{fig:ASGS_Test_case_2}
\end{figure}

\begin{figure}[htbp]
  \centering
  \begin{minipage}[b]{0.46\textwidth}
    \centering
    \includegraphics[width=\textwidth]{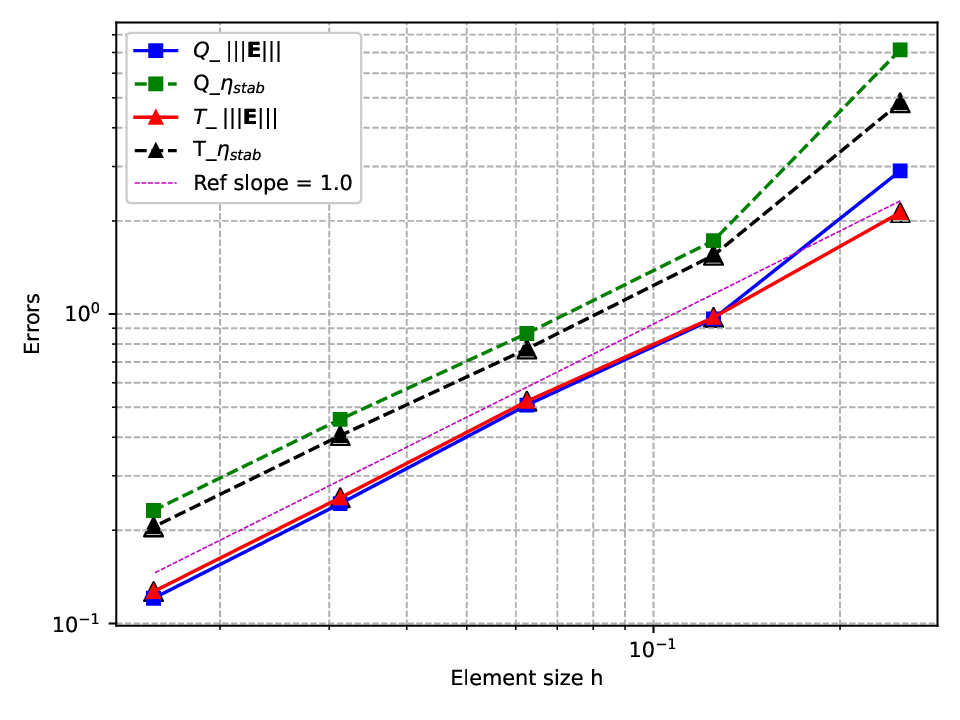}
    \\ \small (a) Stabilized norm ($\eta_{\rm stab}$ vs $|||\bm{E}|||$).
  \end{minipage}%
  \hfill
  \begin{minipage}[b]{0.46\textwidth}
    \centering
    \includegraphics[width=\textwidth]{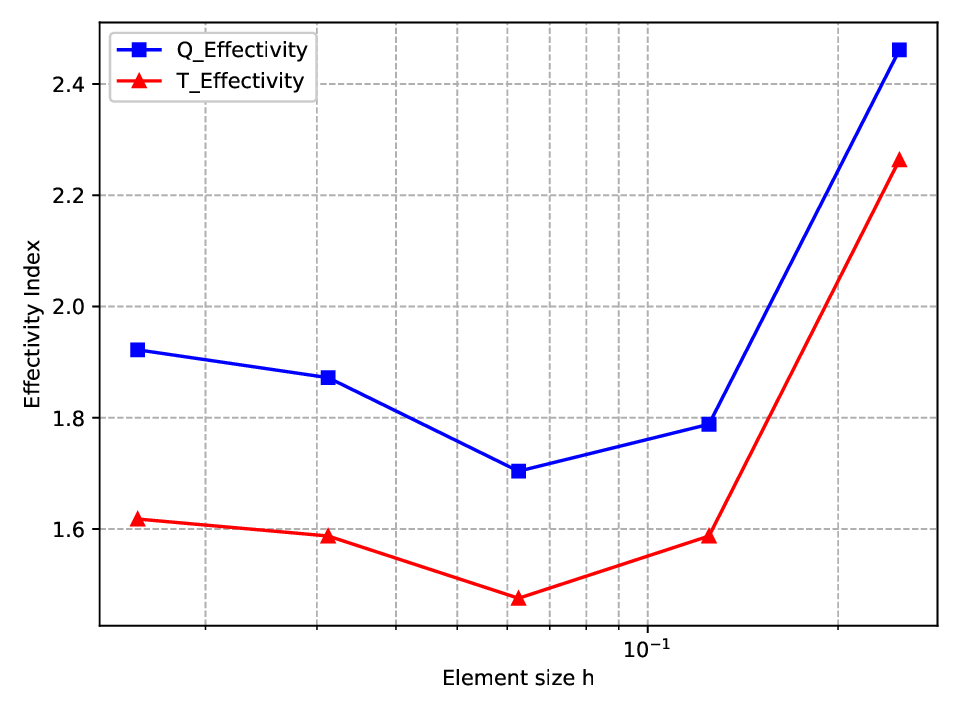}
    \\ \small (b) Effectivity index in stabilized norm.
  \end{minipage}

  \vspace{1ex} 

  \begin{minipage}[b]{0.46\textwidth}
    \centering
    \includegraphics[width=\textwidth]{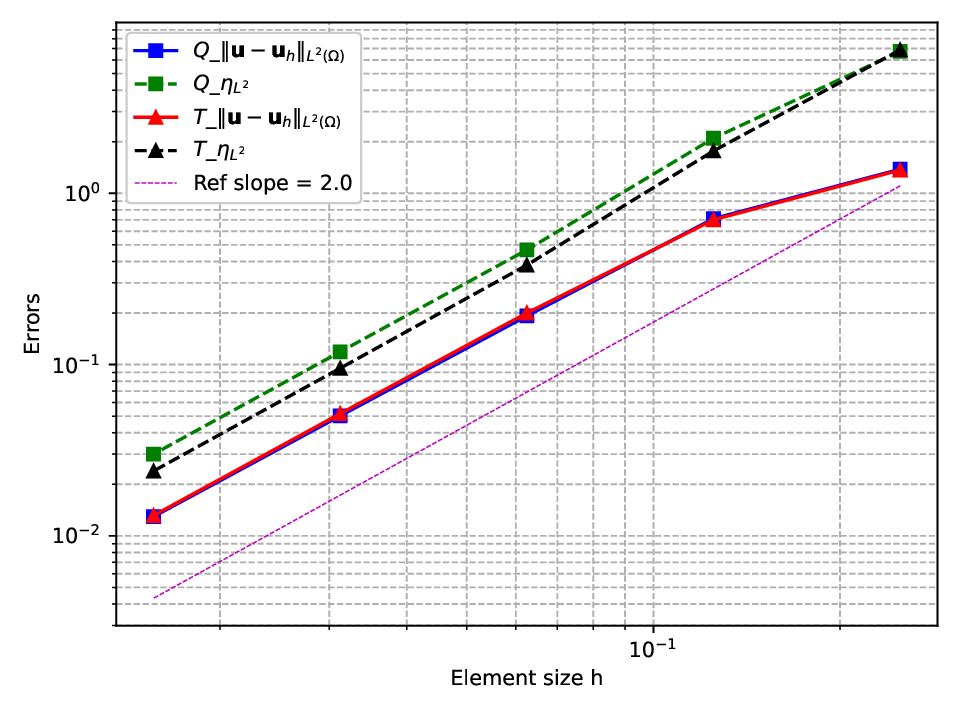}
    \\ \small (c)  $L^2$ norm ($\eta_{L^2}$ vs $E_{L^2}$).
  \end{minipage}%
  \hfill
  \begin{minipage}[b]{0.46\textwidth}
    \centering
    \includegraphics[width=\textwidth]{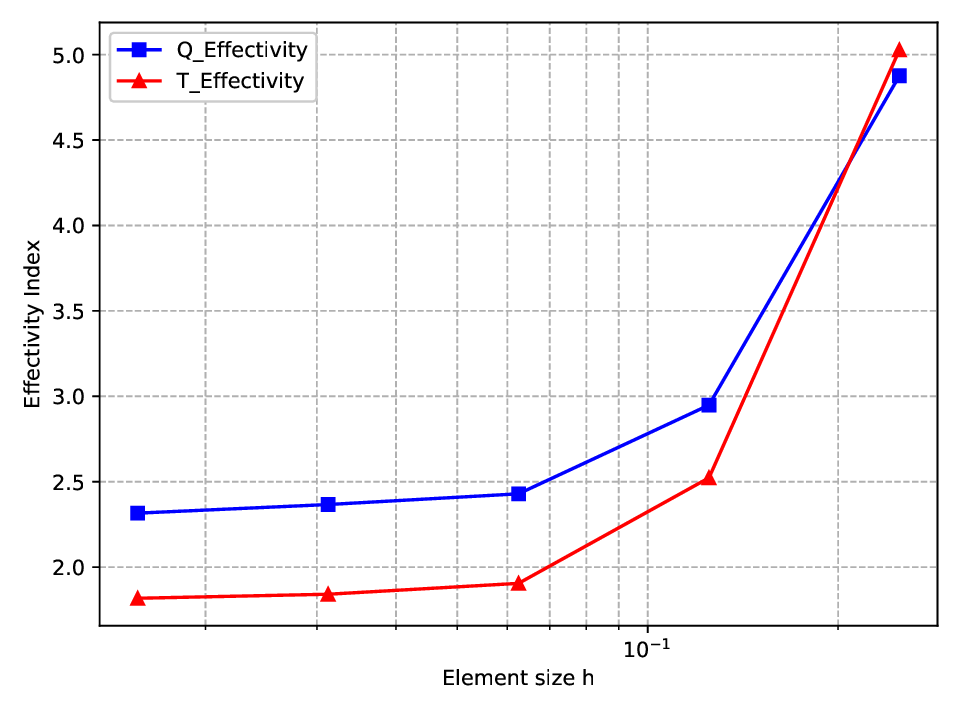}
    \\ \small (d) Effectivity index in $L^2$ norm.
  \end{minipage}

  \begin{minipage}[b]{0.46\textwidth}
    \centering
    \includegraphics[width=\textwidth]{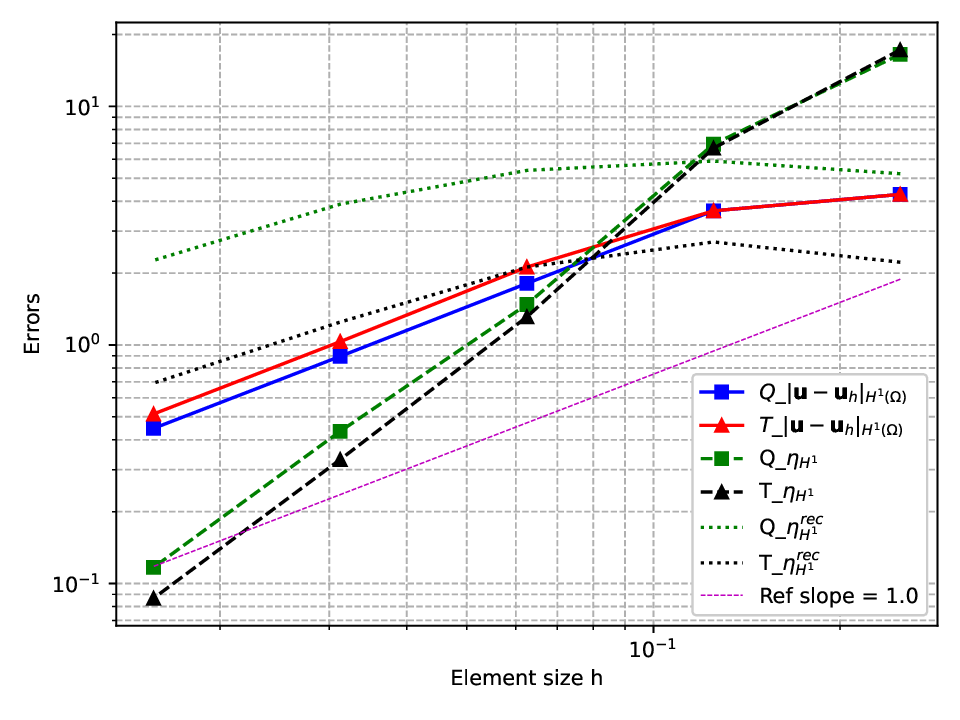}
    \\ \small (e) $H^1$ seminorm ($\eta_{H^1}$, $\eta_{H^1}^{\rm rec}$ vs $E_{H^1}$).
  \end{minipage}%
  \hfill
  \begin{minipage}[b]{0.46\textwidth}
    \centering
    \includegraphics[width=\textwidth]{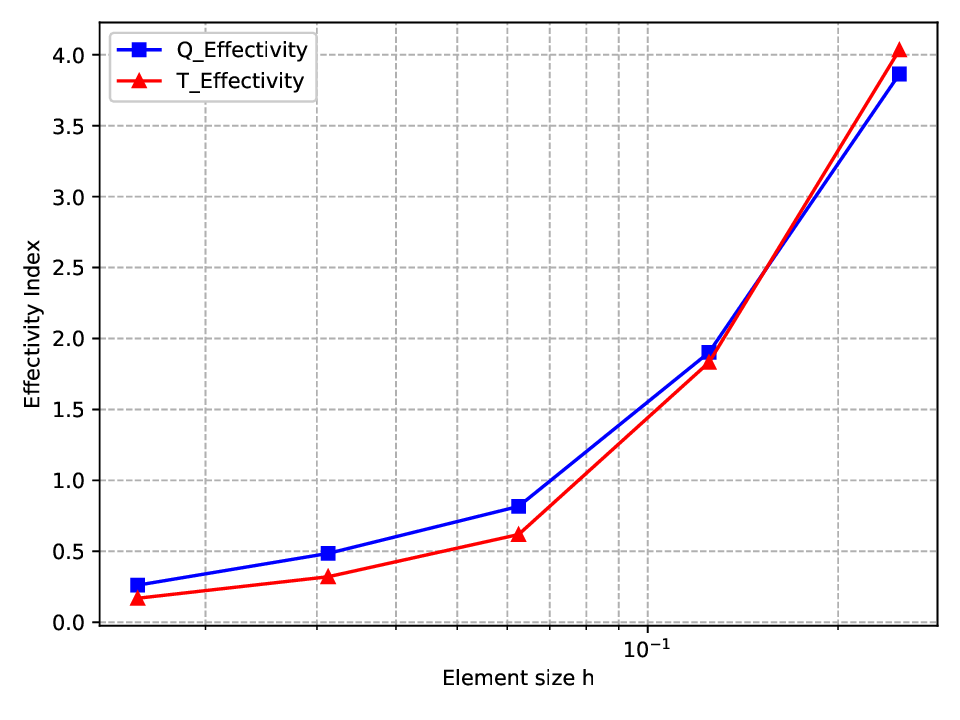}
    \\ \small (f) Effectivity index in $H^1$ seminorm.
  \end{minipage}
  
  \caption{Space convergence (on the left) and effectivity index results (on the right) against the element size of quadrilateral (Q) and triangular (T) linear elements for \textbf{solution 5}  $[\bm{u}_5, p_5]$, $\gamma=18.0$, $\mu  = 1$, and \textbf{ASGS} stabilization technique.}
  \label{fig:ASGS_Test_case_3}
\end{figure}

Now that the effectiveness of the proposed space error estimator alongside the exact error in the stabilized norm has been demonstrated numerically, the goal is to further establish the impact of the projection (OSGS vs ASGS) on the estimator and observe its consistency in the $L^2$ norm and the $H^1$ seminorm to compare with the reference space error estimator in~\cite{Irisarri2021}. Figs.~\ref{fig:OSGS_Test_case_1}, ~\ref{fig:OSGS_Test_case_2}, and ~\ref{fig:OSGS_Test_case_3} present the results of the convergence test using $[\bm{u}_{i}(x,y), p_{i}(x,y)]$,$i=3,4,5$ in different norms for the OSGS stabilization technique, while Figs.~\ref{fig:ASGS_Test_case_1}, ~\ref{fig:ASGS_Test_case_2}, and ~\ref{fig:ASGS_Test_case_3} show the corresponding results for the ASGS technique. The results are computed using linear triangle elements (T) and bilinear quadrilaterals (Q). The viscosity in these cases is $\mu=1$. 

The convergence of the proposed space error estimator in the stabilized norm, shown in Figs.~\ref{fig:OSGS_Test_case_1}-~\ref{fig:ASGS_Test_case_3} (a), closely follows the exact error for quadrilaterals and triangles, for all three manufactured solutions and for both residual projection types, OSGS and ASGS. These confirm that the proposed space error estimator operates reliably and robustly in the stabilized norm. In terms of effectivity index, OSGS (Figs.~\ref{fig:OSGS_Test_case_1}-~\ref{fig:OSGS_Test_case_3} (b)) converges to an effectivity index close to 1, while ASGS converges to values between 1.5 and 2, overestimating the error. This suggests that the orthogonal projection offers a tighter error bound and is more suitable for error estimation. 

It is very important to remark that the proposed scaled norm of the subscales, $\eta_{\mathrm{stab}}$, is an excellent estimator of the error in the stabilized norm, but it is not a suitable estimator for error in other norms. As shown in Figs.~\ref{fig:OSGS_Test_case_1}-~\ref{fig:OSGS_Test_case_3} (c), the $L^2$ norm of the subscales with OSGS has an additional order of convergence compared to the exact error, so it is a poor estimator. On the other hand, it should be noted that the results for the $L^2$ norm of the subscales with ASGS stabilization technique (Fig.~\ref{fig:ASGS_Test_case_1} (c-d)) accurately model the exact error, although there is no formal theory behind the correct behavior of the error estimator in this condition. Regarding the error in the $H^1$ seminorm (Figs.~\ref{fig:OSGS_Test_case_1}-~\ref{fig:ASGS_Test_case_3} (e)), the proposed estimator, $\eta_{H^1}$, has sharper convergence than the real error for OSGS and ASGS, for all manufactured solutions, for triangles and quadrilaterals, making it unsuitable as an error estimator. The reference estimator by Irisarri and Hauke \cite{Irisarri2021} is the correct approach to estimate the error in the $H^1$ seminorm.

\begin{figure}
  \centering
  \includegraphics[width=0.46\textwidth]{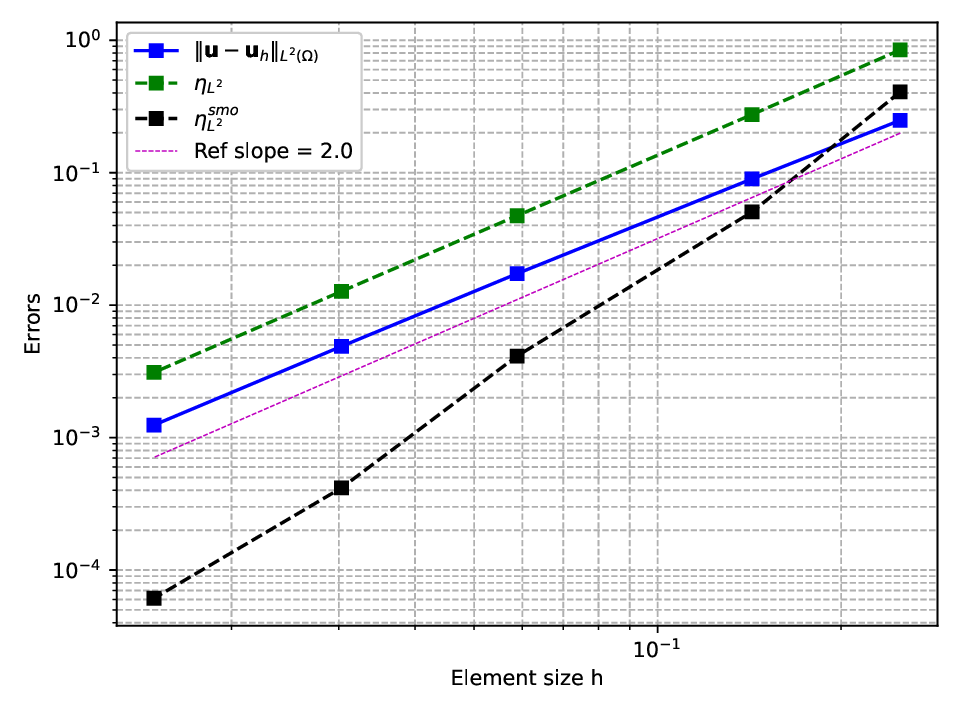}
  \includegraphics[width=0.46\textwidth]{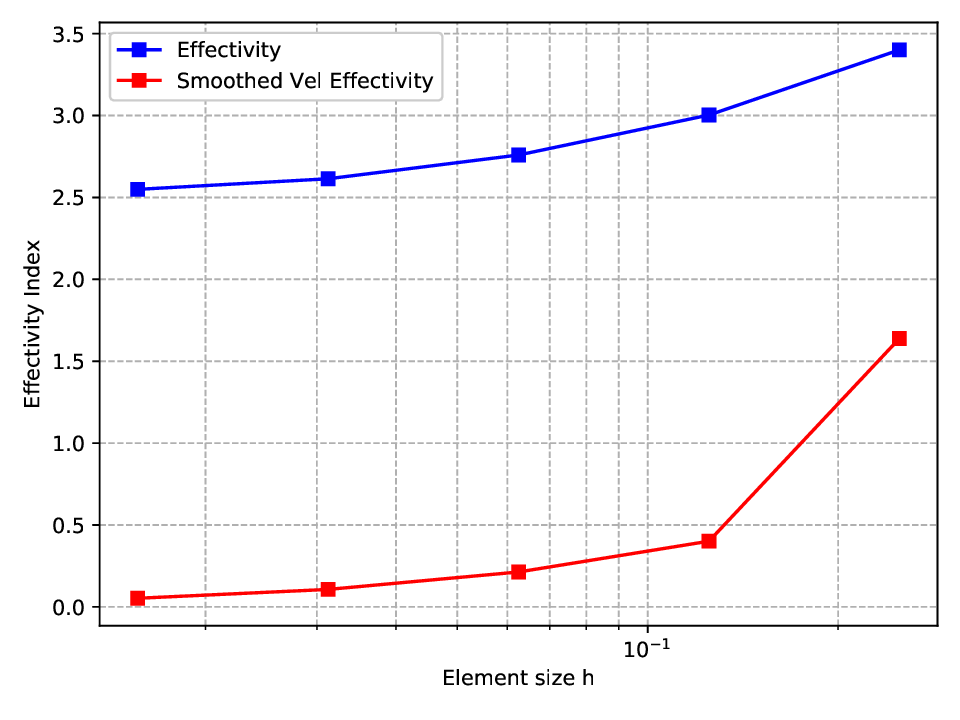}\\

  \caption{Space convergence of errors $\Vert\bm{u}-\bm{u}_{h}\Vert_{L^2(\Omega)}$, $\eta_{L^2}$, and $\eta_{L^2}^{smo}$ (smoothed prior to the norm calculation) on the left and effectivity index on the right $(\mathcal{I}_{\rm eff})$  \textbf{Solution 1}, $[\bm{u}_{1}, p_{1}]$, against the element size of linear quadrilateral elements, for $\gamma=18.0$, $\mu=1$.}
  \label{fig:part3_STST_Space_Error_Indicator_Cases}
\end{figure}

To further analyze the spurious additional convergence, Fig.~\ref{fig:part3_STST_Space_Error_Indicator_Cases} compares the correlation of the proposed space error estimator using ASGS in the $L^2$ norm, where velocity SGSs are calculated normally, compared to the situation where they are smoothed (projected to the nodes) prior to the norm calculation. It can be seen that the smoothing/projection process increases the convergence slope, which is not optimal. This behavior explains the extra convergence slope of the proposed space error estimator in the $H^1$ seminorm, since smoothing the velocity SGSs at the nodes is essential for calculating the seminorm within the linear elements.

%

\subsection{Flow over a backward-facing step at low Reynolds number}
\noindent The first practical problem is the stationary incompressible flow over a backward-facing step which is a benchmark problem in computational fluid mechanics. The problem domain is a two-dimensional backward-facing step with a length of 6 units and a height of 1 unit. The height of the inlet is 0.5, the height of the step is 0.5, and the step is located at 4 units downstream from the inlet \cite{Irisarri2021}. The boundary conditions are as follows: a constant horizontal velocity $ \bm{u} = (1, 0)$ is injected from the inlet boundary, homogeneous Dirichlet boundary conditions are set on the walls, and the right boundary is prescribed as zero traction. Dynamic viscosity coefficients of $\mu = 1$ and $\mu = 0.01$ and a density $\rho = 1$ are used, which results in Reynolds numbers of 1 and 100. The spatial domain is discretized with a uniform mesh of size $h=0.0625$, totaling approximately 1280 linear quadrilateral elements.

\begin{figure}
  \centering
  
  (a)\includegraphics[width=0.46\textwidth]{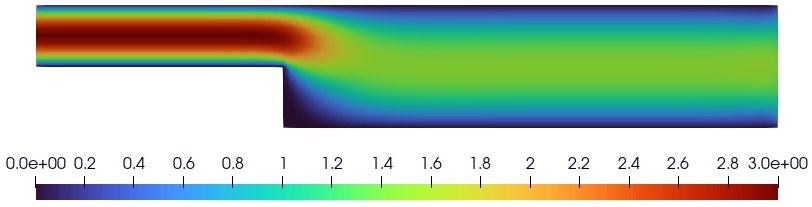}
  (b)\includegraphics[width=0.46\textwidth]{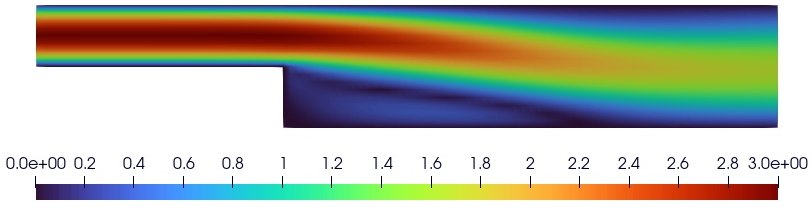}\\
  (c)\includegraphics[width=0.46\textwidth]{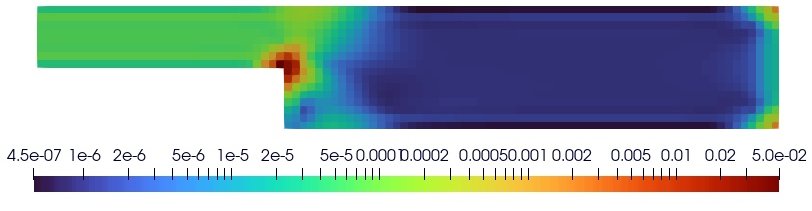}
  (d)\includegraphics[width=0.46\textwidth]{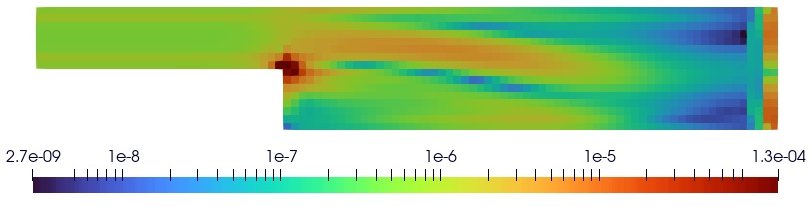}\\
  (e)\includegraphics[width=0.46\textwidth]{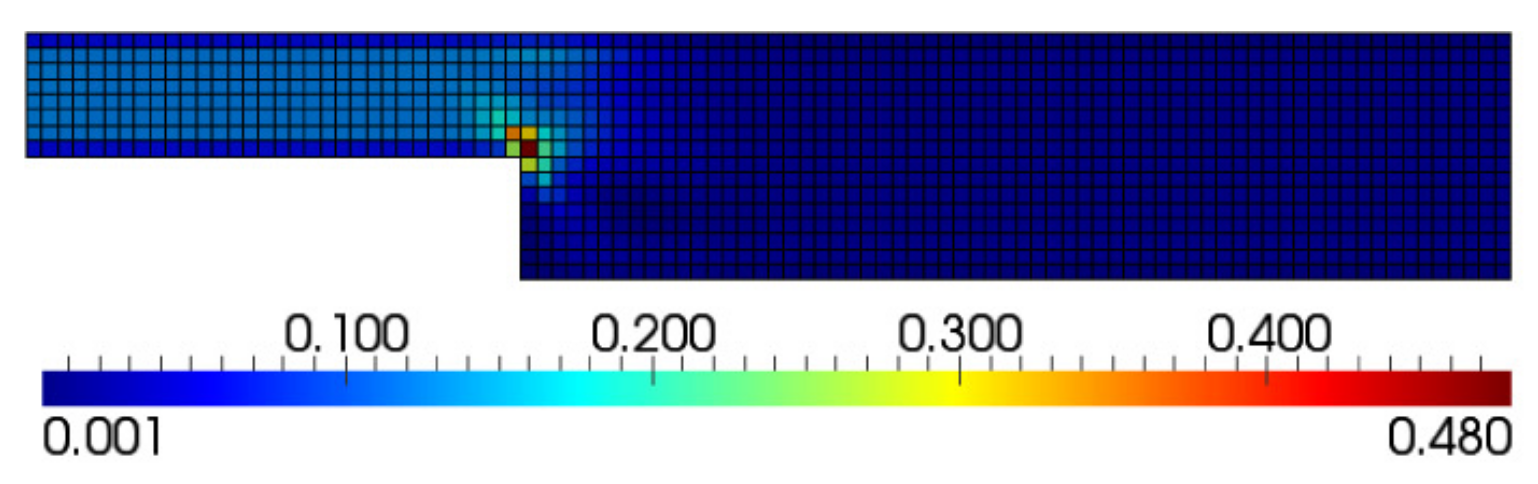}
  (f)\includegraphics[width=0.46\textwidth]{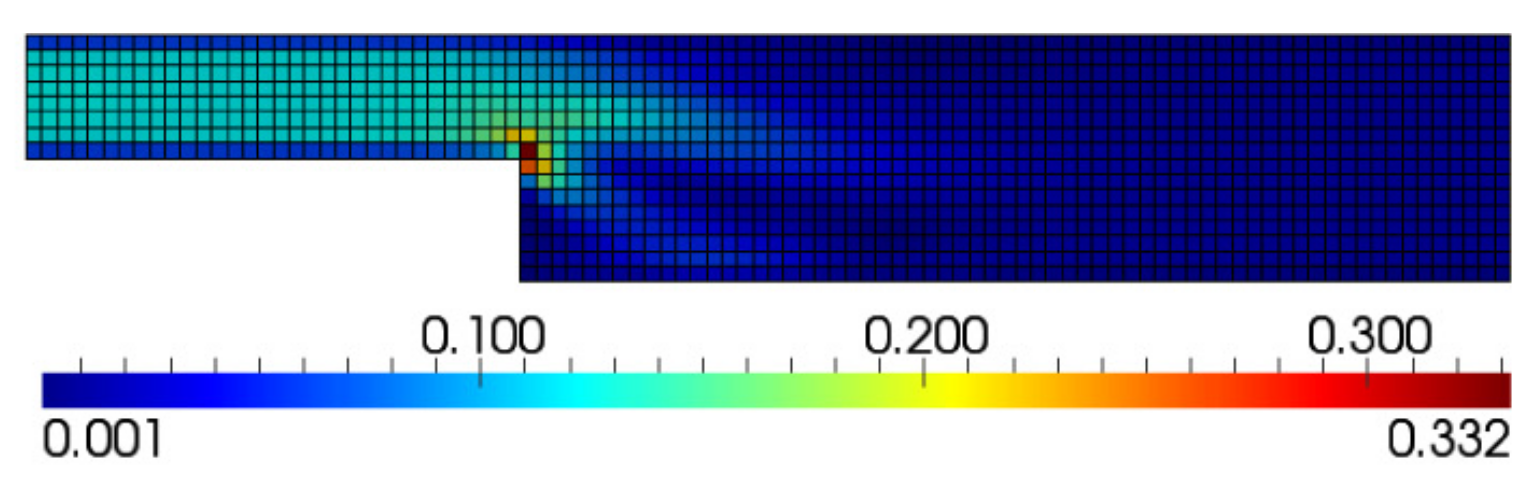}\\
  \caption{Velocity field and error maps for incompressible flow over a backward-facing step at Re = 1 (\textbf{Left}) and Re = 100 (\textbf{Right}): \textbf{(a-b)} velocity magnitudes; \textbf{(c-d)} proposed error estimator, $\eta_{\mathrm{stab}}$; \textbf{(e-f)} Irisarri and Hauke reference error estimator, $\eta_{H^1}^{\rm rec}$.}
  \label{fig:backward-facing_step}
\end{figure}

Fig.~\ref{fig:backward-facing_step} presents the results of the backward-facing step case for Reynolds numbers 1 (on the left) and 100 (on the right), comparing the error maps in the domain obtained with the proposed and reference ~\cite{Irisarri2021} space error estimators. Although these two error estimators are in different norms, it is observed that the error maps follow the same pattern, with the highest error being concentrated at the corner of the step for both Reynolds numbers. The error maps show that both error estimators follow the flow pattern, with errors especially high in the higher-Reynolds case (Figs.~\ref{fig:backward-facing_step}(d) and (f)), where flow recirculation occurs. Moreover, Figs.~\ref{fig:backward-facing_step}(c-d) indicate that the proposed error estimator, $\eta_{\mathrm{stab}}$, is able to detect features like concentrated error at the end of the recirculation region (low velocity but sharp shear layer) and the error arising from the free-traction boundary condition.

\subsection{Flow over a cylinder at Reynolds 500}
\noindent The last case is that of a transient incompressible flow over a cylinder in the transitional regime (Re = 500). The problem domain of this third numerical example is a two-dimensional rectangular domain measuring 16 units in length and 8 units in height, with the bottom-left corner at the coordinate origin. A cylinder with a diameter of 1 unit is located at point (4, 4). For boundary conditions, a constant horizontal velocity  $\bm{u} = (1,0)$ is injected from the left; the vertical component of the velocity is set to zero over both the upper and lower walls (slip condition); zero traction is imposed on the right (outlet) boundary; and for the cylinder, the no-slip condition is prescribed. A dynamic viscosity coefficient $\mu=0.002$ and a density $\rho=1$ are used to ensure a Reynolds number of 500.

\begin{figure}[htbp]
\centering
    (a)\includegraphics[width=0.46\textwidth,height=3.95cm]{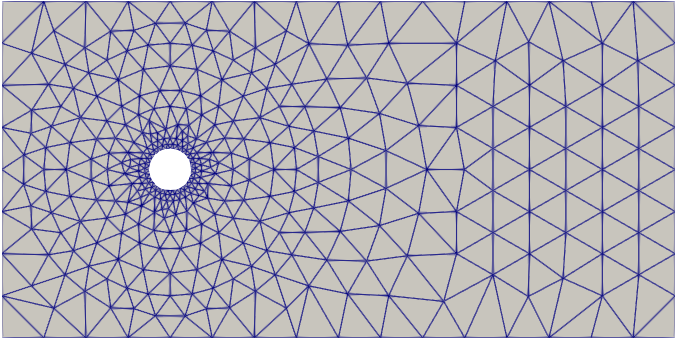}
    (b)\includegraphics[width=0.46\textwidth,height=4.1cm]{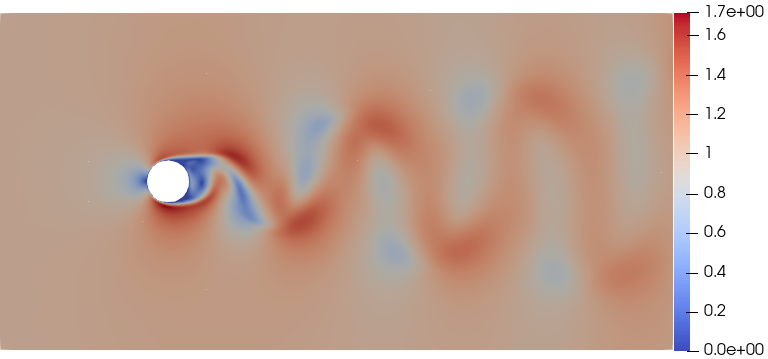} \\
    (c)\includegraphics[width=0.46\textwidth,height=3.95cm]{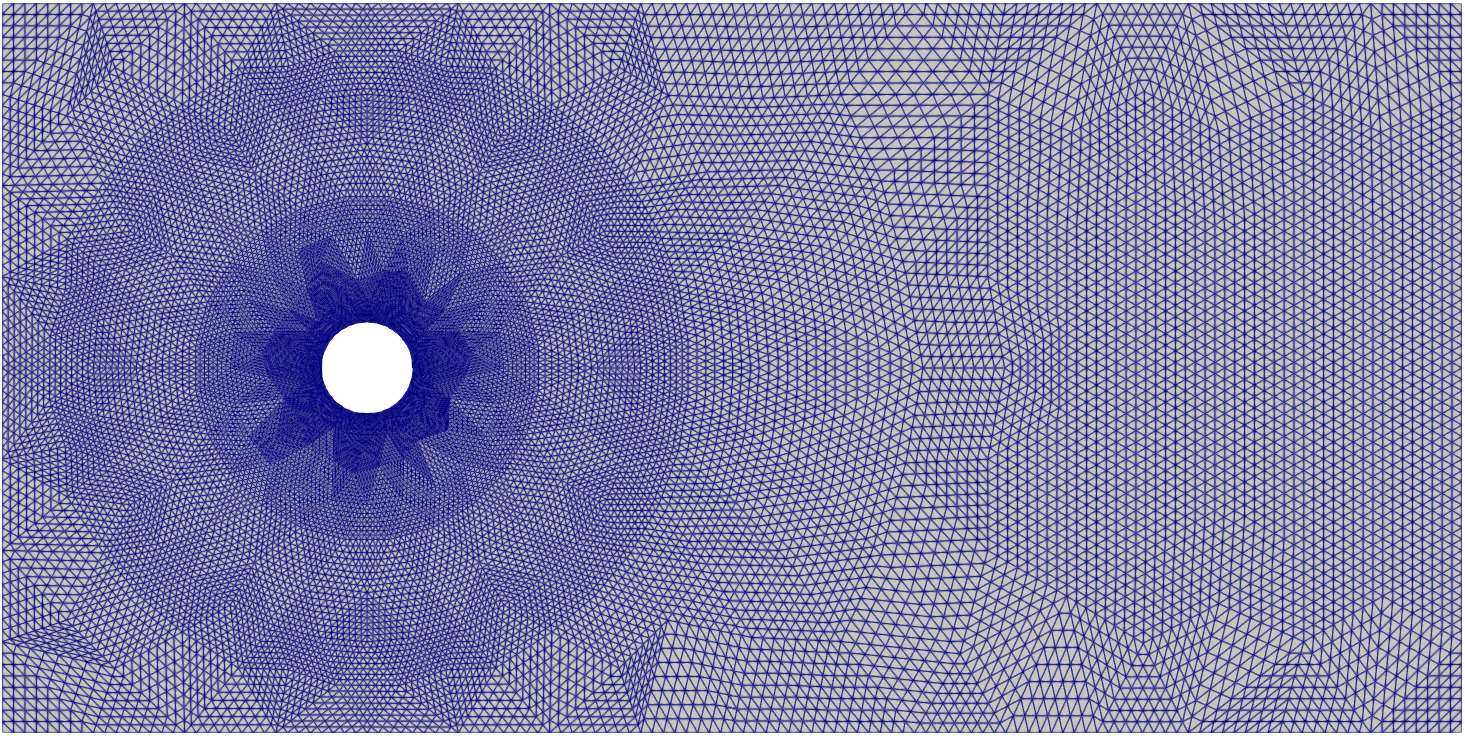} \hspace{0.05cm}
    (d)\includegraphics[width=0.46\textwidth,height=4.1cm]{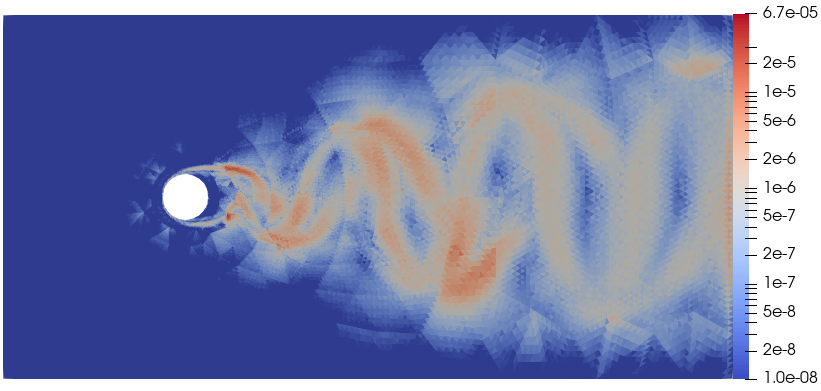}\\
    (e)\includegraphics[width=0.46\textwidth,height=3.95cm]{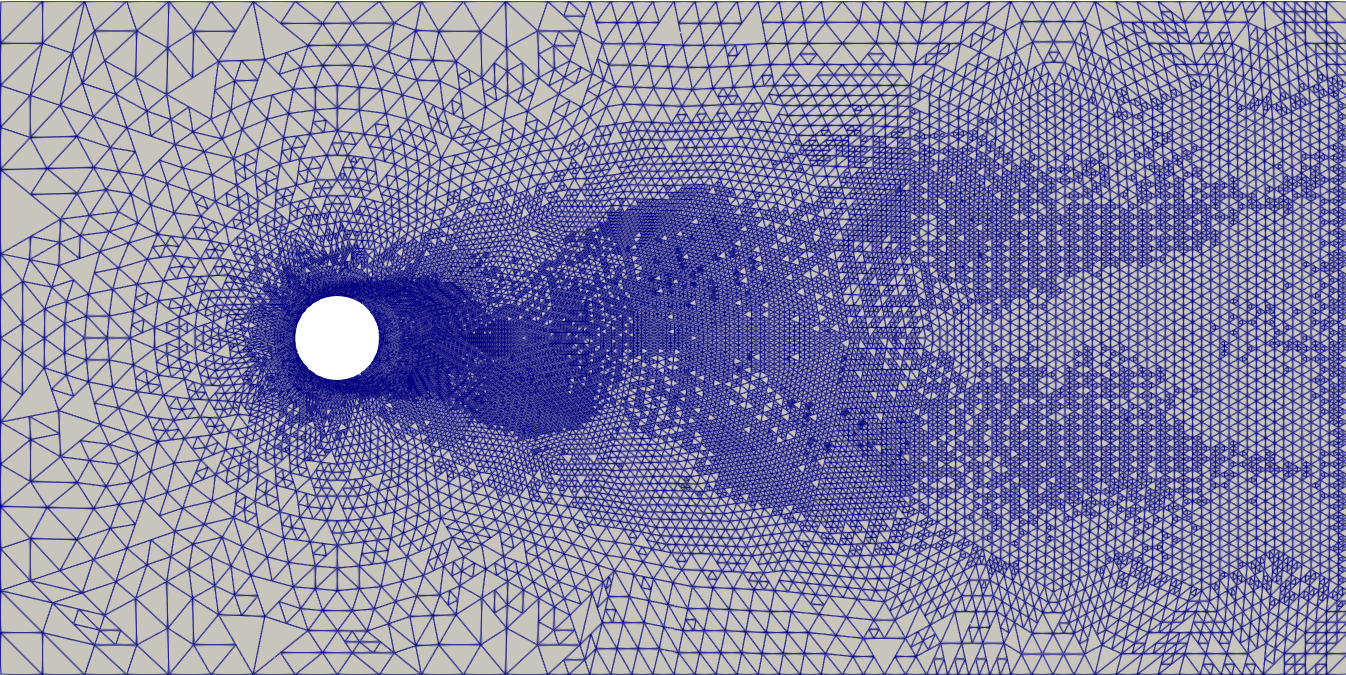} \hspace{0.05cm}
    (f)\includegraphics[width=0.46\textwidth,height=4.1cm]{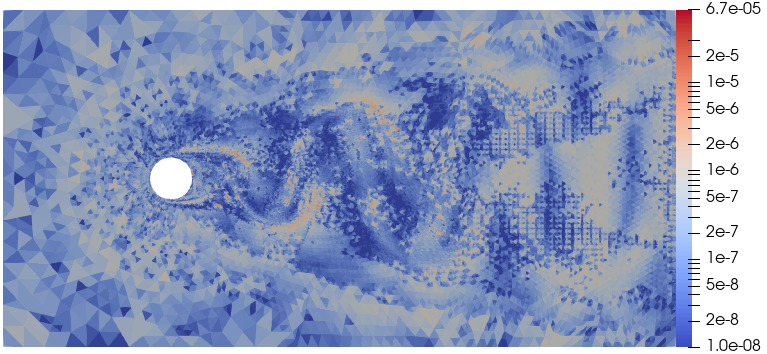}\\
  \caption{Incompressible flow over a cylinder results at time instant $t=97.57$ s: \textbf{(a)} Initial mesh; \textbf{(b)} Velocity magnitude; \textbf{(c)} Mesh after 3 levels of uniform refinement, resulting in 45120 triangular elements; \textbf{(d)} Error estimator, $\eta_{\mathrm{stab}}$, for 3 levels of uniform refinement; \textbf{(e)} Mesh for $\eta_{\mathrm{stab}}$-based refinement case with 45057 triangular elements; \textbf{(f)} Error estimator, $\eta_{\mathrm{stab}}$, for $\eta_{\mathrm{stab}}$-based refinement case.}
  \label{fig:Flow_Over_Cylinder}
\end{figure}

The spatial domain is originally discretized with an unstructured, non-uniformly distributed mesh of 705 linear triangular elements with a characteristic element size of 0.1 around the cylinder, shown in Fig.~\ref{fig:Flow_Over_Cylinder}(a), which will be adaptively refined. 

In this case, the simulation was performed using time discretization with a second-order Backward Difference Formula (BDF2) scheme and a time step size of 0.16678. The time-step size is sufficiently small that the time-discretization error is expected to be small compared to the spatial error for the mesh sizes used here; a formal temporal convergence study is left for future work. Although the methodology of this publication was developed for stationary flow problems, in a computation with sufficiently fine time step, such as this, the spatial error theory remains applicable. The results are presented at a characteristic time instant corresponding to $t=97.57$ s.

To compare the effectiveness of the proposed space error estimator, $\eta_{\mathrm{stab}}$, in a refinement strategy (which is beyond the scope of this paper; the details are therefore omitted), the performance of this method is assessed against the uniform refinement strategy, with both roughly matching the number of elements ($\approx$ 45000 triangular elements) at the time instant when they are compared. The adaptive mesh refinement strategy used is an $h$-refinement hierarchical approach. The computational domain is recursively subdivided by splitting each parent element into child elements (four in 2D for triangles). This structure enables localized adaptive mesh refinement at each time step, guided by the a posteriori error estimator. The hierarchical organization enables efficient multilevel data management, neighbor detection, and refinement operations, and is efficiently implemented and parallelized \cite{baiges2017refficientlib}. In this particular case the acceptable error bound given to the refinement algorithm was $[0.01,0.1]$.

A comparison of Fig.~\ref{fig:Flow_Over_Cylinder} (b) and (d) shows that the space error estimator effectively detects the flow pattern. Uniform refinement is not optimal; instead, regions with vortex shedding, particularly those close to the cylinder where errors are higher, require a more refined mesh. Conversely, the upstream error is almost negligible, 
making a fine mesh in those regions unnecessary. In Fig.~\ref{fig:Flow_Over_Cylinder}(d) and (f), it is clearly shown that in the adaptive refinement case guided by the proposed error estimator, the error is distributed more efficiently and is significantly lower than that of the uniform refinement case. Additionally, the mesh configurations at the time instant for the uniform refinement and algorithmic-refinement cases are displayed in Fig.~\ref{fig:Flow_Over_Cylinder}(c) and (e), respectively.

\section{Conclusion}

In this study, an a posteriori space error estimator grounded in the variational multiscale (VMS) framework
for the stationary incompressible Navier-Stokes equations is introduced, where the modeled SGS in the stabilized norm serves as a reliable error estimator.

First, the proposed space error estimator was evaluated using two different manufactured solutions in both viscous-dominated and convection-dominated flow regimes, with the Reynolds number taken as 0.01, 1.0, and 200. The proposed space error estimator demonstrates robust effectiveness across various conditions, particularly when using the OSGS stabilization technique. It can also detect changes in the flow regime, as indicated by shifts in the convergence slope, an important advantage of this estimator. To further assess its performance, three additional manufactured solutions were employed to compare the proposed space error estimator with a reference one suggested by Irisarri and Hauke~\cite{Irisarri2021}. In all cases, the proposed subscales error estimator successfully represented the exact error in the stabilized norm, as it was designed to do. Furthermore, when comparing the results of OSGS and ASGS, OSGS demonstrates a better effectivity index, as it consistently yields tighter upper and lower bounds than the ASGS approach, which is in agreement with findings in \cite{codina2025posteriorierrorestimatesfinite}.

Regarding other norms, it has been shown that directly taking the subscales used for stabilization as an error estimator is not universal in the $L^2$ norm and is incorrect in the $H^1$ seminorm. It is worth noting that for the ASGS stabilization technique, the proposed error estimator in the $L^2$ norm exhibits convergence behavior consistent with that of the exact error, although there is no theoretical requirement. In the $H^1$ seminorm, it has been confirmed that taking the norm of the gradient of the velocity subscales projected to the nodes is an inadequate estimator of velocity error in $H^1$, regardless of the stabilization technique used. Our study has shown that projecting velocity SGS values onto the nodes, a necessary step for evaluating the $H^1$ seminorm within the linear element, increases the convergence slope. This finding suggests that our SGS approximation is unsuitable for directly estimating the error in the $H^1$ seminorm, making the proposed error estimator unreliable for this purpose. In this sense, the reference estimator by Irisarri and Hauke \cite{Irisarri2021} is an appropriate methodology to estimate velocity error in $H^1$. 

The reference error estimator requires a reconstruction process involving bubble functions and time-averaged constants and provides an accurate estimate of the velocity error in the $H^1$ seminorm. On the other hand, the proposed estimator provides an accurate estimate of the error in the stabilized norm and relies solely on parameters already computed for the stabilization of the governing equations, making it an efficient and practical choice. Both error estimation approaches are valid, but in different norms.

In practical situations such as the flow over a backward-facing step, both error estimators identify similar regions of high error and can therefore guide adaptive refinement in comparable ways, despite estimating different norms. The principal advantage of the proposed error estimator is that it bypasses the reconstruction process and relies solely on parameters already available from the stabilization of the solver. Another advantage is that the proposed estimator, because it is designed for the stabilized norm, represents the global error of the problem and not only the error on a single unknown, such as the velocity.

Furthermore, the ultimate objective is to employ the proposed spatial error indicator within an adaptive refinement algorithm. The flow over a cylinder problem serves as a case study to illustrate the effectiveness of the proposed error estimator in guiding adaptive mesh refinement. As shown in Fig. ~\ref{fig:Flow_Over_Cylinder}, the results demonstrate that the error indicator facilitates a uniform distribution of errors within user-defined bounds, thereby validating its practical applicability in adaptive algorithms.

\section*{Acknowledgements}

I. Mart\'inez-Su\'arez gratefully acknowledges the support received from the Ag\`encia de Gesti\'o d’Ajuts i de Recerca (AGAUR) through the predoctoral FI grant 2022-FI-B-00517. R. Codina gratefully acknowledges the support received through the ICREA Aca\-dèmia Research Program of the Catalan Government.

\section*{Declarations}

\subsection*{Conflict of interest}

The authors declare that they have no known competing financial interests or personal relationships that could have appeared to influence the work reported in this paper.

\subsection*{Replication of results}

All the information required for replicating the paper results was duly presented. Data files for the results are available upon request to the authors.

\bibliographystyle{ieeetr}
\bibliography{references.bib}
\end{document}